\documentclass[journal]{IEEEtran}

\usepackage{color}
\usepackage[pdftex]{graphicx}
\usepackage{algorithm,algorithmic}
\usepackage{array}
\usepackage{amssymb,amsmath,color} 
\usepackage{subfig}
\usepackage{epstopdf}
\usepackage{multirow}

\newcommand{\real}{{\mathbb{R}}}
\newcommand{\subj}{\text{subj. to}}

\newcommand{\ponewt}{\texttt{PRONTO}}
\newcommand{\force}{F}

\newtheorem{remark}{Remark}

\begin{document}

\title{Minimum-time trajectory generation \\ for quadrotors in constrained environments}

\author{Sara~Spedicato,
  and~Giuseppe~Notarstefano,~\IEEEmembership{Member,~IEEE}
  \thanks{S. Spedicato and G. Notarstefano are with the Department of
    Engineering, Universit\`a del Salento, Via per Monteroni, 73100 Lecce,
    Italy, e-mail:
    \texttt{name.lastname@unisalento.it.}}
  \thanks{This result is part of a project that has received funding from the
    European Research Council (ERC) under the European Union's Horizon 2020
    research and innovation programme (grant agreement No 638992 - OPT4SMART).}
\thanks{A short, preliminary version of this work was presented at the IEEE Conference on Decision and Control 2016. Differences with that work include: (i) a more comprehensive treatment of the proposed strategy, (ii) an extended version of the strategy involving the computation of collision free regions shaped by obstacles, (iii) numerical computations for more general scenarios and (iv) an experimental test on a nano-quadrotor.}}

\maketitle


\begin{abstract}
  In this paper, we present a novel strategy to compute minimum-time
  trajectories for quadrotors in constrained environments. In particular, we
  consider the motion in a given flying region with obstacles and take into
  account the physical limitations of the vehicle.  Instead of approaching the
  optimization problem in its standard time-parameterized formulation, the
  proposed strategy is based on an appealing re-formulation. Transverse
  coordinates, expressing the distance from a frame path, are used to
  parameterize the vehicle position and a spatial parameter is used as
  independent variable. This re-formulation allows us to (i) obtain a fixed
  horizon problem and (ii) easily formulate (fairly complex)
  position constraints. The effectiveness of the proposed strategy is
  proven by numerical computations on two different illustrative
  scenarios. Moreover, the optimal trajectory generated in the second scenario
  is experimentally executed with a real nano-quadrotor in order to show its
  feasibility.
\end{abstract}



\IEEEpeerreviewmaketitle


\section{Introduction}
\IEEEPARstart{N}{}umerous applications involving Unmanned Aerial Vehicles (UAVs), and in
particular quadrotors, require them to move inside areas characterized by
physical boundaries, obstacles and even tight space constraints (as e.g., urban
environments) in order to accomplish their robotics tasks.  Such applications
are, for example, structural inspections, transportation tasks, surveillance and
search and rescue missions.  Trajectory generation, a core step for physical
task realization \cite{dadkhah2012survey}, becomes extremely challenging in this
scenario.  A physically realizable trajectory must satisfy (i) the (nonlinear)
system dynamics, (ii) the physical limits of the vehicle, such as the maximum
thrust, and (iii) the position constraints.
Although safety (ensured by a feasible trajectory) is the primary requirement
for all applications, trajectory optimization 
is becoming necessary in different application domanis.
The cost to minimize can be, for example, the time to execute a
maneuver (in a search and rescue scenario), the energy consumption (during long
endurance missions), or the ``distance'' from a desired unfeasible state-input curve
(during inspections).  The further requirement of performance optimization
poses an additional challenge in the trajectory generation problem.
		
The problem of computing optimal paths (or trajectories) 
for UAVs (e.g., \cite{bottasso2008path} and \cite{ambrosino2009path}) has received significant attention and a number of
algorithms for quadrotors have been proposed to accomplish complex tasks, e.g., landing on a
moving target \cite{herisse2012landing} and blind navigation in unknown
populated environments \cite{naldi2015robust}.  Focusing on collision avoidance,
two different approaches, namely reacting or planning, can be applied.  The
reactive approach is based on navigation laws preventing from possible
collisions. It can be performed, e.g., modulating the velocity reference
\cite{hou2016dynamic}, selecting ad-hoc reference way-points
\cite{furci2015plan} and defining an harmonic potential field
\cite{masoud2015plan}.  On the contrary, the planned approach deals with a
problem involving dynamics and state-input constraints with (possibly) a
performance criterion to optimize. The majority of the planning algorithms
regarding quadrotors, such as \cite{cowling2010direct},
\cite{mellinger2011minimum}, \cite{bouktir2008trajectory}, \cite{van2013time}, \cite{chen2016online},
takes advantage of the differential flatness property and relies on approximations via motion primitives. 
When dealing with obstacle dense environments, trajectory generation is often
performed using a decoupled approach (\cite{bry2015aggressive},
\cite{koyuncu2008probabilistic}, 
\cite{bouffard2009hybrid}).  In
a first stage, a collision-free path is generated by sampling-based path planning
algorithms, such as the Rapidly-exploring Random Tree (RRT) in
\cite{bry2015aggressive, koyuncu2008probabilistic} or the Probabilistic Roadmap
(PRM) in 
\cite{bouffard2009hybrid}, and without the dynamics constraint.  In a
second stage, an optimal trajectory (satisfying the system dynamics) is
generated from the collision-free path. Optimization techniques such as
\cite{cowling2010direct}, \cite{mellinger2011minimum},
\cite{bouktir2008trajectory} can be used at this stage. 
In order to overcome the limitations due to the
decoupled approach, 
a variant of the RRT algorithm is developed in \cite{devaurs2015optimal}, 
an approximated dynamics with an a-posteriori correction is used in \cite{allen2016real}
and a space-parameterized problem reformulation, suitable for modeling complex flight scenarios, is adopted in \cite{van2013time}. 
Differently from the previous works, in \cite{hehn2015real} the structure of the minimum-time trajectories is found by the Pontryagin's minimum principle. Nevertheless, position constraints are not considered.
Finally, in \cite{augugliaro2012generation} a discretized simple point-mass dynamics and
approximated convex constraints are considered. 
The approximation of non-convex constraints into convex ones is also used in \cite{augugliaro2013dance}, in which a sequential convex programming approach is used to achieve a collision free motion for dancing quadrotors.
	
Our main contribution is the design of an optimization framework to generate
feasible minimum-time quadrotor trajectories in structured environments as,
e.g., rooms, corridors, passages, or urban areas.  Our strategy computes
optimal trajectories that satisfy the quadrotor nonlinear dynamics. 
The strategy can be applied to general models, which may be more complicated than 
the differentially flat ones.
  Instead of addressing the minimum-time problem in its standard free-horizon
  formulation, we derive a fixed-horizon reformulation in which transverse
  coordinates, expressing the ``transverse" distance from a frame path, are used
  to parameterize the vehicle position. The resulting problem, having a spatial
  parameter as independent variable, is easier to solve than the
  time-parametrized one.  Position constraints can be easily added into the
  reformulated problem by defining the constraint boundaries as a function of
  the spatial parameter and shaping them according to the presence of obstacles.
  Approximate solutions to the infinite-dimensional optimization problem are
  numerically computed by combining the Projection Operator Newton method for
Trajectory Optimization (PRONTO) \cite{hauser2002projection} with a barrier
function approach \cite{hauser2006barrier}. This method generates
  trajectories in a numerically stable manner and guarantees recursive
  feasibility during the algorithm evolution, i.e., at each algorithm iteration
  a system trajectory is available. Moreover the approximated solution always
  satisfies the constraints since the barrier function approach is an interior
  function method.
As an additional contribution, we present numerical computations to show the
effectiveness of the proposed strategy on two challenging scenarios. In the
first one, the moving space is delimited by rooms with obstacles of different shapes.
In the second scenario, the constrained
environment is a tubular region delimited by hula hoops.
The optimal minimum-time trajectory related to this second scenario is experimentally
performed on our nano-quadrotor testbed.

Our algorithm compares to the literature in the following way.
The majority of works, such as \cite{cowling2010direct, mellinger2011minimum, bouktir2008trajectory, van2013time, chen2016online}, uses the differential
flatness to avoid the integration of nonlinear differential equations, to reduce
the order of the problem and to simplify the definition of constraints
\cite{cowling2010direct}.  
On the contrary, our strategy does not rely on the differential flatness hypothesis and thus it can be applied to more complex models.
In the previously cited works, the optimization problem is posed in the flat output space, where outputs are approximated using motion primitives, such as
polynomial functions \cite{cowling2010direct, mellinger2011minimum,
chen2016online}, B-splines \cite{bouktir2008trajectory}, or ``convex
combinations of feasible paths" \cite{van2013time}. The optimization variables are thus the parameters of the motion primitives.
Differently from these works, we do not rely on motion primitives: the state-input trajectory is the optimization variable in our problem formulation. Similarly to the problem formulation in \cite{hauser2006motorcycle}, our reformulated minimum-time problem has a spatial parameter, instead of time, as independent variable. While in \cite{hauser2006motorcycle} the maximum velocity profile (for a given path) is computed for a motorcycle model by using a quasi-static approximation of the dynamics, we optimize the whole state-input trajectory and we consider the full nonlinear dynamics of the quadrotor.
Finally, other optimization strategies using the PRONTO method are \cite{hauslerenergy} and \cite{rucco2015virtual}, which aim to compute respectively minimum-energy trajectories for two-wheeled mobile robots and minimum-distance trajectories (from an unfeasible desired maneuver) for UAVs.
Differently from these works, we consider a more general three-dimensional space with position constraints and we
reformulate the minimum-time problem by using the transverse coordinates.

  The paper is organized as follows. In Section \ref{sec:problem_formulation} we
  present the standard formulation of the optimization problem we aim to
  solve. In Section \ref{sec:strategy} our trajectory generation strategy, based
  on an appealing reformulation of the problem, is described.  Finally, in
  Section \ref{sec:computations}, we provide numerical computations and
  experiments, and discuss interesting features of the computed minimum-time
  trajectories.


\section{The quadrotor minimum-time problem}
\label{sec:problem_formulation}
We first briefly introduce the quadrotor model used in the paper and then recall
the standard problem formulation.
	
\subsection{Quadrotor model}	
The quadrotor dynamics can be described by the so called vectored-thrust
dynamical model in \cite{hua2013introduction}, where the gravity is the only external force and the generated torque does not influence the translational dynamics,
i.e., 
\begin{align}
\dot{\boldsymbol{p}} &= {\text{\textbf{v}}}\label{eq:state_pos}\\
\boldsymbol{\dot{\text{\textbf{v}}}}& = g \boldsymbol{e}_3 - \frac{\force}{m} R(\boldsymbol{\Phi})  \boldsymbol{e}_3 \label{eq:state_vel}\\
\dot{\boldsymbol{\Phi}}&=J(\boldsymbol{\Phi}) \boldsymbol{\omega} \label{eq:state_ang}\\
\boldsymbol{\dot{\omega}} &= -I^{-1}\hat{\boldsymbol{\omega}} I \boldsymbol{\omega} + I^{-1} \boldsymbol{\gamma}. \label{eq:state_omega}
\end{align}
with $\boldsymbol{p}=[p_1 \; p_2 \: p_3]^T$, $\text{\textbf{v}}=[\text{v}_1 \; \text{v}_{2} \; \text{v}_3]^T$, $\boldsymbol{\Phi} = [\varphi \; \theta \; \psi]^T$, where $\varphi$, $\theta$, $\psi$ are respectively the roll, pitch and yaw angles, and $\boldsymbol{\omega} = [p \; q \; r]^T$.
The symbols in equations (\ref{eq:state_pos}-\ref{eq:state_omega}) are defined in the following table, where $\mathcal{F}_i$ and $\mathcal{F}_b$ respectively denote the inertial and the body frame.
\begin{table}[hb]
\begin{center}
  \caption{Nomenclature}\label{tb:symbols}
\begin{tabular}{cc}
\hline\\[-2ex]
$\boldsymbol{p} \in \real^3$ & position vector expressed in $\mathcal{F}_i$ \\[0.5ex]
$\text{\textbf{v}}\in \real^3$ & velocity vector expressed in $\mathcal{F}_i$\\[0.5ex]
$\boldsymbol{\Phi} \in \real^3$ & vector of angles (yaw-pitch-roll w.r.t. current frame)\\[0.5ex]
$R(\boldsymbol{\Phi}) \!\in\! SO(3)$ & rotation matrix to map vectors in $\mathcal{F}_b$ into vectors in $\mathcal{F}_i$\\[0.5ex]
$\boldsymbol{\omega} \in \real^3$ & angular rate vector expressed in $\mathcal{F}_b$\\[0.5ex]
$\hat{\boldsymbol{\omega}} \in so(3)$ & skew-symmetric matrix associated to $\boldsymbol{\omega}$ \\[0.5ex]
$J(\boldsymbol{\Phi}) \in \mathbb{R}^{3 \times 3}$ & matrix mapping $\boldsymbol{\omega}$ into $\dot{\boldsymbol{\Phi}}$ \\[0.5ex]
$m \in \real$ & vehicle mass\\[0.5ex]
$I \in \mathbb{R}^{3 \times 3}$ & inertia matrix\\[0.5ex]
$g \in \real$ & gravity constant\\[0.5ex]
$\boldsymbol{e}_3\in \real^3$ & vector defined as $\boldsymbol{e}_3:=[ 0 \; 0 \; 1]^T$ \\[0.5ex]
$\force \in \real$ & thrust\\[0.5ex]
$\boldsymbol{\gamma} \in \real^3$ & torque vector\\[0.5ex]
\hline
\end{tabular}
\end{center}
\end{table}
	
For the vehicle maneuvering, we adopt a cascade control scheme with an off-board
position/attitude control loop and an on-board angular rate controller.
Assuming that the \emph{virtual} control input $\boldsymbol{\omega}$ is tracked by the
on-board angular rate controller, we restrict our trajectory generation problem
on the position/attitude subsystem (\ref{eq:state_pos}-\ref{eq:state_ang}),
which can be written in state-space form as
\begin{equation}
  \dot{\boldsymbol{x}}(t) = f(\boldsymbol{x}(t),\boldsymbol{u}(t)),
  \label{eq:state_space}
\end{equation}
with state $\boldsymbol{x} = [\boldsymbol{p}^T \; \text{\textbf{v}}^T \; \boldsymbol{\Phi}^T]^T,$ input $\boldsymbol{u} = [\boldsymbol{\omega}^T \; \force]^T$ and suitably defined $f$. 
	
\subsection{Quadrotor minimum-time problem: standard formulation}
We deal with the following optimal control problem:
\begin{align}
  \begin{split}
    \min_{\boldsymbol{x}(\cdot),\boldsymbol{u}(\cdot),T} &\;\;  T\\
    \subj &\;\; \dot{\boldsymbol{x}}(t) = f(\boldsymbol{x}(t),\boldsymbol{u}(t)), \quad \boldsymbol{x}(0)= \boldsymbol{x_0} \; \text{\emph{(dynamics)}}\\
    & \;\; \boldsymbol{x}(T) \in {X}_T \; \text{\emph{(final constraint)}}\\
    & \;\; |p(t)| \leq p_{max} \; \text{\emph{(roll rate)}}\\
    & \;\; |q(t)| \leq q_{max} \; \text{\emph{(pitch rate)}}\\
    & \;\; |r(t)| \leq r_{max} \; \text{\emph{(yaw rate)}}\\
    & \;\;  0 < \force_{min} \leq \force(t) \leq \force_{max} \; \text{\emph{(thrust)}} \\
    & \;\; |\varphi(t)| \leq \varphi_{max}(t) \; \text{\emph{(roll angle)}}\\
    & \;\; |\theta(t)| \leq \theta_{max}(t) \; \text{\emph{(pitch angle)}}\\
    & \;\; |\psi(t)| \leq \psi_{max}(t) \; \text{\emph{(yaw angle)}}\\
    & \;\; c_{obs}(\boldsymbol{p}(t)) \leq 0 \; \text{\emph{(position constraints)}},
\end{split}
  \label{eq:mintime_standard}
      \end{align}
where ${X}_T \subset \real^9$ is a desired final region, $p_{max}$, $q_{max}$ and $r_{max}$ are bounds on roll, pitch and yaw rate,
respectively, $\force_{min}$ and $\force_{max}$ are lower and upper bounds on thrust,
$\varphi_{max}(\cdot)$, $\theta_{max}(\cdot)$ and $\psi_{max}(\cdot)$ are bounds on roll-pitch-yaw angles, and
$c_{obs} : \real^3 \rightarrow \real$ represents position
constraints.
The $t$-dependent constraints in \eqref{eq:mintime_standard} hold for all $t \in [0,T]$, with the exception of $c_{obs}(\boldsymbol{p}(t)) \leq 0$, which holds for all $t \in [0,T)$.
The bounds on the angular rates avoid fast solutions. The vehicle thrust is also limited:
quadrotor vehicles can only generate positive thrust and the maximum rotor speed
is limited. Furthermore, constraints on roll and pitch angles are imposed into the optimization
problem in order to avoid acrobatic vehicle configurations and to satisfy $\theta \neq \pm \frac{\pi}{2}$ (which makes the matrix $J(\boldsymbol{\Phi})$ in \eqref{eq:state_ang} always well defined).
Time dependent boundaries can be used for roll and pitch constraints when the vehicle has to move through small passages.
The constraint on the yaw angle may be useful in applicative scenarios in which a sensor, e.g., a camera, is provided onboard the vehicle and needs to be pointed toward a target region. The position constraints take into account physical boundaries (possibly shaped by the presence of obstacles) and may also represent GPS denied areas or spaces with limited communication.

\section{Minimum-Time Trajectory Generation Strategy}
\label{sec:strategy}

In this section, we describe our strategy to compute minimum-time
trajectories.

Minimum-time problem \eqref{eq:mintime_standard} is difficult to solve, since it
is a constrained, free-horizon problem (time $T$ is an optimization
variable). For this reason, instead of directly designing an algorithm to solve
problem \eqref{eq:mintime_standard}, we provide a strategy to obtain an
equivalent, but computationally more appealing, fixed-horizon formulation.
In the following, we give an informal idea of the strategy steps to derive the new problem formulation.  First, we
define a \emph{frame path} as a (purely geometric) curve in $\real^3$ used to
express the quadrotor position in terms of new coordinates. That is, as depicted in Figure \ref{fig:manreg}, the
position is identified by the arc-length of the point on the path at minimum
distance and by two transverse coordinates expressing how far the quadrotor
position is from the curve.
Second, we rewrite the dynamics in terms of the transverse coordinates and show
that it depends on time only through the arc-length time-evolution. Thus, by
using the arc-length as independent variable, we obtain a ``space-dependent''
transverse dynamics. 
Third, the time $T$ can be expressed itself as a function of the
arc-length over a fixed ``spatial'' horizon $[0,L]$, with $L$ being the total
length of the frame path. Thus, minimizing $T$ can be rewritten as minimizing
an integral function over the fixed spatial interval $[0,L]$. Similarly,
pointwise constraints can be written in terms of the transverse coordinates and
as function of the arc-length. 

The resulting fixed-horizon optimal control problem is solved by using
the Projection Operator Newton method (PRONTO), \cite{hauser2002projection},
combined with a barrier function approach to handle the constraints,
\cite{hauser2006barrier}.

We provide a detailed and formal explanation of the strategy steps in the
following subsections.
	
\subsection{Frame path}
\label{sec:3a}	
The first step of the strategy is the 
generation
of an arc-length parameterized frame path $\bar{\boldsymbol{p}}_f(s)$,
$\forall {s} \in [0,L]$, where $s$ is the arc-length of the path and $L$ is
its total length.  In the following, we denote the arc-length parameterized
functions with a bar, and the derivative with respect to the arc-length with a
prime, i.e.,
$\bar{\boldsymbol{p}}_f'(s) :=
d\bar{\boldsymbol{p}}_f(s)/ds$.
The frame path $\bar{\boldsymbol{p}}_{f}(\cdot)$ has to be locally a
non-intersecting $C^2$ curve with non-vanishing
$\bar{\boldsymbol{p}}'_{f} (\cdot)$.  Note that the frame path is only a
geometric path and it is not required to satisfy the position constraints. A possibility is the computation of the frame path as a $C^{\infty}$ geometric curve, e.g., using arctangent
functions as in our numerical computations.  More details on the frame path used
for our numerical computations will be given in Section \ref{sec:computations}.

The frame path is used to parameterize the inertial position of the vehicle in the new transverse coordinates, as will be clear later. In order to define the transverse coordinates, we consider the Serret-Frenet frame, whose origin has $\bar{\boldsymbol{p}}_{f}(s)$ as coordinates, and defined $\forall s \in [0,L]$. 
In particular, the tangent, normal and bi-normal vectors, respectively $\bar{\boldsymbol{t}}(s), \bar{\boldsymbol{n}}(s),
\bar{\boldsymbol{b}}(s)$, are defined, with components in the inertial frame, as
\begin{align}
  \bar{\boldsymbol{t}}(s)&:=\bar{\boldsymbol{p}}'_f(s),\\
  \bar{\boldsymbol{n}}(s)&:=\frac{\bar{\boldsymbol{p}}''_f(s)}{\bar{k}(s)},\\
  \bar{\boldsymbol{b}}(s)&:= \bar{\boldsymbol{t}}(s) \times
                           \bar{\boldsymbol{n}}(s),
\end{align}
where $\bar{k}(s):=\Vert \bar{\boldsymbol{p}}''_f(s)\Vert_2$ is the
curvature of $\bar{\boldsymbol{p}}_f(\cdot)$ at $s$. 
Moreover, we define the rotation matrix 
\begin{align}
\bar{R}_{SF}:=[\: \bar{\boldsymbol{t}} \; \bar{\boldsymbol{n}} \; \bar{\boldsymbol{b}}\:]
\label{eq:Rsf}
\end{align}
mapping
vectors with components in the Serret-Frenet frame into vectors with
components in the inertial frame. 
According to the Serret-Frenet formulas \cite{PSM:08}, the arc-length derivative of the Serret-Frenet rotation matrix is 
\begin{equation}
  \bar{R}'_{SF}(s) =
  \bar{R}_{SF}(s)
  \left[
    \begin{array}{ccc}
      0 & -\bar{k}(s)  & 0 \\
      \bar{k}(s) & 0 & -\bar{\tau}(s) \\
      0 & \bar{\tau}(s) & 0 \\
    \end{array}
  \right],
  \label{eq:Rsf'*dot_s_lqr}
\end{equation}
where $\bar{\tau}(s):=\bar{\boldsymbol{n}}(s) \; \bar{\boldsymbol{b}}'(s)$ is
the torsion of $\bar{\boldsymbol{p}}_f(\cdot)$ at $s$.

\subsection{Transverse dynamics}
\label{sec:transv}
The second step of the strategy is the derivation of the transverse dynamics by using the transverse coordinates defined with respect to the frame path $\bar{\boldsymbol{p}}_f(\cdot)$. In order to rewrite the standard dynamics (\ref{eq:state_pos}-\ref{eq:state_ang}) into the transverse dynamics, we proceed as follows.

\begin{figure}[htbp]
  \begin{center}
    \includegraphics[scale=1.3]{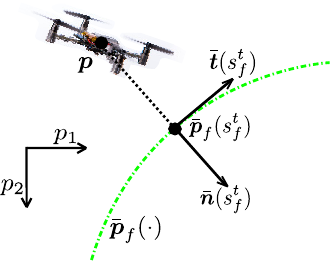}\\
      \caption{Selection of the arc-length $s$ identifying the point on the frame path at minimum distance from the quadrotor position at the time instant $t$.}
    \label{fig:manreg}
  \end{center}
\end{figure}

First, we design a change of coordinates from the inertial position $\boldsymbol{p} \in \real^3$ to
the transverse coordinate vector $\boldsymbol{w} \in \real^2$, such that 
$\boldsymbol{w} = [w_1 \; w_2]^T$, 
where $w_1$ and $w_2$ are the transverse coordinates.
Let us consider the quadrotor center of mass with position $\boldsymbol{p}(t)$.
As depicted in Figure \ref{fig:manreg}, its orthogonal projection on the frame path identifies a point with position $\bar{\boldsymbol{p}}_f(s_f(t))$,
where the function $s_f : \real_{0}^+ \rightarrow \real_{0}^+$ is
\begin{align}
s_f(t):=\text{arg} \min_{s \in \real_0^+}
\|\boldsymbol{p}(t)-\bar{\boldsymbol{p}}_f(s)\|^2.
\label{eq:pi-lpqr}
\end{align}
For simplicity, in the following we use $s_f^t := s_f(t)$ and $\dot{s}_f^t := \dot{s}_f(t)$.
Note that, the minimizing arc-length is unique provided that $\bar{\boldsymbol{p}}_f(\cdot)$
is locally a non-intersecting $C^2$ curve with non-vanishing $\bar{\boldsymbol{p}}'_f (\cdot)$.
By mapping 
$\boldsymbol{p}-\bar{\boldsymbol{p}}_f(s_f^t)$ into a vector with components in the Serret-Frenet frame attached to $\bar{\boldsymbol{p}}_f(s_f^t)$, we obtain 
\begin{align}
\boldsymbol{d} &:= \bar{R}_{SF}(s_f^t)^T (\boldsymbol{p}-\bar{\boldsymbol{p}}_f(s_f^t)).
\label{eq:d}
\end{align}  
Noticing that the component related to the tangent vector is always zero by construction, we define the components $w_1$ and $w_2$
of the transverse error vector $\boldsymbol{w}$
as, respectively, the second and third components of $\boldsymbol{d}$, i.e.,
\begin{align}
\begin{split}
  w_1& := \bar{\boldsymbol{n}}(s_f^t)^T (\boldsymbol{p}-\bar{\boldsymbol{p}}_f(s_f^t)),\\
  w_2& := \bar{\boldsymbol{b}}(s_f^t)^T
                          (\boldsymbol{p}-\bar{\boldsymbol{p}}_f(s_f^t)),
\end{split}
\label{eq:x_change_coordinates}
\end{align}
and thus obtaining 
\begin{align}
\boldsymbol{d} = [0 \; w_1 \; w_2]^T.
\label{eq:d-def}
\end{align}
\vspace{-0.2cm}

Second, we rewrite equation \eqref{eq:state_pos} using $\boldsymbol{w}$ instead of $\boldsymbol{p}$.
We note that, the invertible function $s_f(\cdot)$ provides a change of variables from the time $t$ to the arc-length $s$.
A generic arc-length function $\bar{\alpha}(\cdot)$ can be expressed as the time function $\bar{\alpha}(s_f(\cdot))$ and its time derivative is $\frac{d\bar{\alpha}(s_f(t))}{dt} = \bar{\alpha}'(s_f^t) \; \dot{s}_f^t$. 
Let us rewrite equation \eqref{eq:state_pos}.
By using 
equation \eqref{eq:d}, the position of the quadrotor center of mass $\boldsymbol{p}(t)$, at time
instant $t$, can be written as
\begin{align}
  \boldsymbol{p}(t) = \bar{\boldsymbol{p}}_f(s_f^t)+\bar{R}_{SF}(s_f^t) \; \boldsymbol{d}(t).
  \label{eq:y=y_xi+Rsf*d_lqr}
\end{align}
Differentiating \eqref{eq:y=y_xi+Rsf*d_lqr} with respect to time, since \eqref{eq:state_pos} holds, we get
\begin{equation}
  \text{\textbf{v}}(t)=\bar{\boldsymbol{p}}'_f(s_f^t) \; \dot{s}_f^t +
  {\bar{R}}_{SF}'(s_f^t) \; \dot{s}_f^t \; \boldsymbol{d}(t) +
  \bar{R}_{SF}(s_f^t) \; \boldsymbol{\dot{d}}(t).
  \label{eq:dot_y=y_xi+Rsf*d_lqr}
\end{equation}
Multiplying both sides of equation \eqref{eq:dot_y=y_xi+Rsf*d_lqr} by $\bar{R}^{T}_{SF}$, 
using \eqref{eq:Rsf'*dot_s_lqr}, \eqref{eq:d-def} and $\bar{\boldsymbol{p}}'_f(s_f^t) = \bar{R}_{SF}(s_f^t)
[1 \; 0 \; 0]^T$,
we get
\begin{equation*}
\left[
    \begin{array}{c}
0 \\
\dot{w}_1(t)\\
\dot{w}_2(t)
    \end{array}
\right]
+
\dot{s}_f^t
\left[
\begin{array}{ccc}
      1-\bar{k}(s_f^t) w_1(t)\\
      -\bar{\tau}(s_f^t) w_2(t)\\
      \bar{\tau}(s_f^t) w_1(t)\\
\end{array}
\right]    
-
\bar{R}_{SF}^T(s_f^t) \text{\textbf{v}}(t)
=
0,
\label{eq:dot_y2}
\end{equation*}
i.e., using \eqref{eq:Rsf},
\begin{align}
\dot{s}_{f}^t &= \frac{\bar{\boldsymbol{t}}(s_{f}^t)^T \text{\textbf{v}}(t)}{1-\bar{k}(s_{f}^t) w_1(t)} \label{eq:s_}\\
  \dot{w}_1(t) &= \bar{\boldsymbol{n}}(s_{f}^t)^T \text{\textbf{v}}(t) +\bar{\tau}(s_{f}^t) \dot{s}_{f}^t w_2(t) \label{eq:w1_}\\
  \dot{w}_2(t) &= \bar{\boldsymbol{b}}(s_{f}^t)^T \text{\textbf{v}}(t) -\bar{\tau}(s_{f}^t) \dot{s}_{f}^t w_1(t). \label{eq:w2_}
\end{align}


Third and final, we rewrite equations \eqref{eq:w1_}, \eqref{eq:w2_}, \eqref{eq:state_vel}, \eqref{eq:state_ang}, 
by using the arc-length $s$ as independent variable.
Let us denote by $\bar{t}_f : \real_{0}^+ \mapsto \real_{0}^+$ the inverse function of $s_f : \real_{0}^+ \mapsto \real_{0}^+$, satisfying $t = \bar{t}_f(s_f^t)$.
Due to the invertibility of $s_f(\cdot)$, a generic time function $\alpha(\cdot)$ can be expressed as the arc-length function $\alpha(\bar{t}_f(\cdot))$ and, defining $\bar{\alpha} := \alpha \circ \bar{t}_f$, 
we have $\alpha(t) = \bar{\alpha}(s_f^t)$. 
In particular,
\begin{align}
\boldsymbol{w}(t)&=\bar{\boldsymbol{w}}(s_f^t), \quad \text{\textbf{{v}}}(t)=\bar{\text{\textbf{{v}}}}(s_f^t), \quad \boldsymbol{\Phi}(t)=\bar{\boldsymbol{\Phi}}(s_f^t), \label{eq:bar-w}\\
\boldsymbol{\omega}(t)&=\bar{\boldsymbol{\omega}}(s_f^t), \quad \force(t)=\bar{\force}(s_f^t). \label{eq:bar-omega}
\end{align}
Deriving with respect to time equations \eqref{eq:bar-w}, we get
\begin{align*}
\begin{split}
\dot{\boldsymbol{w}}(t)&=\bar{\boldsymbol{w}}'(s_f^t) \dot{s}_f^t,
\quad 
\dot{\text{\textbf{{v}}}}(t) = \bar{\text{\textbf{{v}}}}'(s_f^t) \dot{s}_f^t, 
\quad
\dot{\boldsymbol{\Phi}}(t) = \bar{\boldsymbol{\Phi}}'(s_f^t) \dot{s}_f^t,
\end{split}
\end{align*}
and equations \eqref{eq:w1_},\eqref{eq:w2_},\eqref{eq:state_vel},\eqref{eq:state_ang} become
\begin{align}
\begin{split}
\label{eq:transv-all-v1}
\bar{w}'_1(s_f^t) &= \bar{\boldsymbol{n}}(s_f^t)^T {\text{\textbf{v}}}(t) \; \frac{1}{\dot{s}_f^t} +\bar{\tau}(s_f^t) {w}_2(t),\\
\bar{w}'_2(s_f^t) &= \bar{\boldsymbol{b}}(s_f^t)^T {\text{\textbf{v}}}(t) \;\frac{1}{\dot{s}_f^t} -\bar{\tau}(s_f^t) {w}_1(t),\\
\bar{\text{\textbf{v}}}'(s_f^t) &= (g \boldsymbol{e}_3 - \frac{{\force}(t)}{m} R({\boldsymbol{\Phi}}(t)) \boldsymbol{e}_3)  \;\frac{1}{\dot{s}_f^t},\\
\bar{\boldsymbol{\Phi}}'(s_f^t) &=J ({\boldsymbol{\Phi}}(t)) {\boldsymbol{\omega}(t)} \; \frac{1}{\dot{s}_f^t}.
\end{split}
\end{align}
Using \eqref{eq:s_}, \eqref{eq:bar-w} and \eqref{eq:bar-omega},
equations \eqref{eq:transv-all-v1} depend on time only through the variable $s_f^t$.
Thus, we can rewrite the dynamics in the arc-length, $s \in [0,L]$, domain. Formally, considering $s$ as
the independent variable, we get the \emph{transverse dynamics}
\begin{align}
\begin{split}
\label{eq:transv-all}
  \bar{w}'_1 &= \bar{\boldsymbol{n}}^T \bar{\text{\textbf{v}}} \; \frac{1-\bar{k} \bar{w}_1}{\bar{\boldsymbol{t}}^T \bar{\text{\textbf{v}}}} +\bar{\tau} \bar{w}_2,\\
  \bar{w}'_2 &= \bar{\boldsymbol{b}}^T \bar{\text{\textbf{v}}} \; \frac{1-\bar{k} \bar{w}_1}{\bar{\boldsymbol{t}}^T \bar{\text{\textbf{v}}}} -\bar{\tau} \bar{w}_1,\\
\bar{\text{\textbf{v}}}'
  &= 
(g \boldsymbol{e}_3 - \frac{\bar{\force}}{m} R(\bar{\boldsymbol{\Phi}}) \boldsymbol{e}_3)  \; \frac{1-\bar{k} \bar{w}_1}{\bar{\boldsymbol{t}}^T \bar{\text{\textbf{v}}}},\\
  \bar{\boldsymbol{\Phi}}' &=J (\bar{\boldsymbol{\Phi}}) \bar{\boldsymbol{\omega}} \; \frac{1-\bar{k} \bar{w}_1}{\bar{\boldsymbol{t}}^T \bar{\text{\textbf{v}}}}.
  \end{split}
\end{align}
Note that the dependence by $s$ is omitted for simplicity.
Equations \eqref{eq:transv-all} can be written in state-space form as
\begin{equation}
\bar{\boldsymbol{x}}'_w (s) = \bar{f}(\bar{\boldsymbol{x}}_w(s),\bar{\boldsymbol{u}}(s)),
  \label{eq:tran_dynamics}
\end{equation}
with state $\bar{\boldsymbol{x}}_w = [\bar{\boldsymbol{w}}^T \; \bar{\text{\textbf{v}}}^T \; \bar{\boldsymbol{\Phi}}^T]^T$, input $\bar{\boldsymbol{u}} = [\bar{\boldsymbol{\omega}}^T \; \bar{\force}]^T$ {and suitable $\bar{f}$}.

\begin{remark}
The general theory regarding the transverse coordinates is introduced in
\cite{AB-JH:94} and used to design a maneuver regulation controller for a
bi-dimensional case in \cite{AS-JH-AB:13}. Differently from \cite{AS-JH-AB:13},
we use the transverse coordinates in a more general three-dimensional case and
in order to develop a trajectory optimization strategy rather than a
controller.
\end{remark}

\subsection{Arc-length parameterization of cost and constraints}
\label{sec:cost-constraints}
The third step of the strategy consists into the reformulation of cost and constraints in problem \eqref{eq:mintime_standard} by using the new (arc-length dependent) variables $\bar{\boldsymbol{x}}_w$ and $\bar{\boldsymbol{u}}$.

The cost functional in \eqref{eq:mintime_standard}, i.e., $T = \int_{0}^T 1 \; dt$,
is rewritten into an arc-length parameterization by considering the change of variable from $t$ to $s$, i.e., 
$$
\int_{0}^T 1 \; dt = \int_{s_f(0)}^{s_f(T)} \bar{t}_f^{\;'}(s) \; ds.
$$
Since $\frac{d\bar{t}_f(s_f(t))}{dt} = \bar{t}_f^{\;'}(s_f^t) \dot{s}_f^t$
and
$\bar{t}_f(s_f^t) = t$, 
we get
\begin{align}
\bar{t}_f^{\;'}(s_f^t) = 1/\dot{s}_f^t
\label{eq:t_f_pr}
\end{align}
with $\dot{s}_f^t$ as in \eqref{eq:s_}.
Since $w_1(t) = \bar{w}_1(s_f^t)$ and $\text{\textbf{v}}(t) = \bar{\text{\textbf{v}}}(s_f^t)$, as in \eqref{eq:bar-w}, equation \eqref{eq:t_f_pr} can be written as 
\begin{align}
  \bar{t}_f^{\;'}(s_f^t) &= \frac{1-\bar{k}(s_{f}^t) \bar{w}_1(s_{f}^t)}{\bar{\boldsymbol{t}}(s_f^t)^T \bar{\text{\textbf{v}}}(s_f^t)}, \label{eq:int_}
\end{align} 
where all the variables depend on time only through $s_f^t$.
Thus, we can rewrite \eqref{eq:int_} in the arc-length, $s \in [0,L]$, domain, obtaining
\begin{align}
  \bar{t}_f^{\;'}(s) &= \frac{1-\bar{k}(s) \bar{w}_1(s)}{\bar{\boldsymbol{t}}(s)^T \bar{\text{\textbf{v}}}(s)}. 
\label{eq:int_s}
\end{align}
Finally, since $s_f(0) = 0$, $s_f(T) = L$, and \eqref{eq:int_s} holds, we rewrite the cost functional in \eqref{eq:mintime_standard} as
\begin{equation}
\int_0^L \!\!\! \quad \frac{1-\bar{k}(s) \bar{w}_1(s)}{\bar{\boldsymbol{t}}(s)^T \bar{\text{\textbf{v}}}(s)} \; ds.
\label{eq:cost_s}
\end{equation}
Notice that, according to \eqref{eq:cost_s}, the hypothesis $\bar{\boldsymbol{t}}(s)^T \bar{\text{\textbf{v}}}(s) \neq 0$, has to be satisfied $\forall s \in [0,L]$, i.e., the velocity projected on the tangent vector of the frame path has to be not null.

\begin{figure}[htbp]
	\begin{center}
	{\includegraphics[scale=0.5]{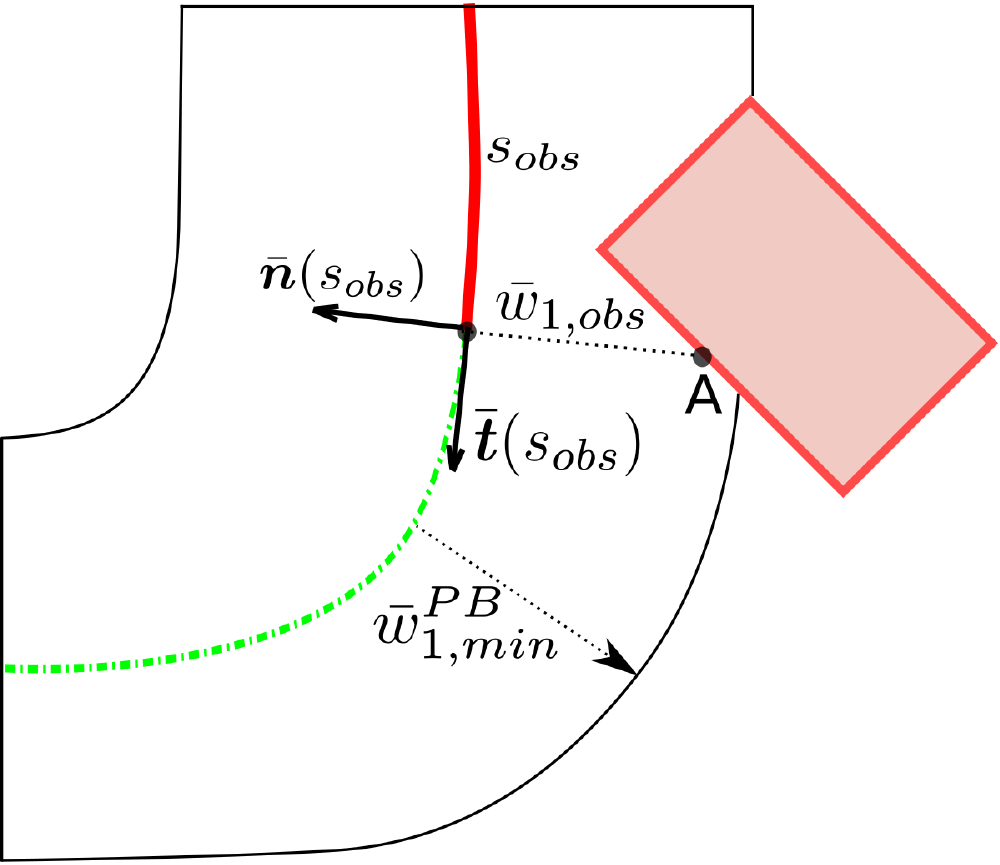}\vspace{0.1cm}\label{fig:path_s2}}\\
	\caption{Representation of $\bar{w}_{1,obs}$ related to a point A on the obstacle boundary. The frame path is depicted in red (portion identifying $s_{obs}$) and dot-dashed green.}
	\label{fig:w1_obs}
	\end{center}	
\end{figure}

The constraints in \eqref{eq:mintime_standard} are rewritten into an arc-length parameterization suitable to apply the barrier function approach \cite{hauser2006barrier}.
The constraint $\boldsymbol{x}(T) \in \boldsymbol{X}_T$ is written in the form 
\begin{align}
c_f(\bar{\boldsymbol{x}}_w(L)) \leq 0,
\label{eq:cT}
\end{align}
with scalar components
\begin{align}
  c_{f,i}(\bar{x}_{w_i}(L)) = \Big( \frac{2 \; \bar{x}_{w_i}(L) - (\bar{x}_{w_i,max}+\bar{x}_{w_i, min})}{(\bar{x}_{w_i, max}-\bar{x}_{w_i, min})} \Big)^2 -1, 
  \label{eq:c_fi}
\end{align}
$\forall i= 1,...,8$, where 
$\bar{x}_{w_i}$ is the $i$-th component of $\bar{\boldsymbol{x}}_{w}$,
and $\bar{x}_{w_i, min}$ and $\bar{x}_{w_i, max}$ are the bounds on the final states.
The constraints on the angular rates, thrust and roll-pitch-yaw angles are rewritten by using 
equations \eqref{eq:bar-w}, \eqref{eq:bar-omega} and reparameterizing the time-dependent bounds 
$\varphi_{max}(t)$, $\theta_{max}(t)$ and $\psi_{max}(t), \; \forall t \in [0,T]$, by the arc-length $s$. 
Thus, we have
\begin{align}
\begin{split}
\label{eq:constr_s}
&\Big( \frac{\bar{p}(s)}{p_{max}} \Big)^2 \!\!\!-1 \leq 0, \;\;\;  
\Big( \frac{\bar{q}(s)}{q_{max}} \Big)^2 \!\!\!-1 \leq 0, \;\;\; 
\Big( \frac{\bar{r}(s)}{r_{max}} \Big)^2 \!\!\!-1 \leq 0, \\
&\Big( \frac{\bar{\varphi}(s)}{\bar{\varphi}_{max}(s)} \Big)^2\!\!\!\!-\!1\!\leq\!0,
\Big( \frac{\bar{\theta}(s)}{\bar{\theta}_{max}(s)} \Big)^2\!\!\!\!-\!1\!\leq\!0,
\Big( \frac{\bar{\psi}(s)}{\bar{\psi}_{max}(s)} \Big)^2\!\!\!\!-\!1\!\leq\!0,\\
& \Big( \frac{2 \bar{\force}(s) - (\force_{max}+\force_{min})}{(\force_{max}-\force_{min})}\Big)^2 \!\!\!-1 \leq 0. 
\end{split}
\end{align}
As regards the position constraints $c_{obs}(\boldsymbol{p}(t)) \leq 0$, they are written in the generic form
\begin{align}
\bar{c}_{obs}(\bar{w}_1(s),\bar{w}_2(s)) \leq 0,
\label{eq:pos_constr}
\end{align}
which can be particularized according to the shape of the flying region.
For environments with circular sections, the inequality \eqref{eq:pos_constr} becomes
\begin{align}
\Big( \frac{\sqrt{\bar{w}_1^2(s)+\bar{w}_2^2(s))}}{\bar{r}_{obs}(s)} \Big)^2 -1 \leq 0,
\label{eq:circ_constr}
\end{align}
where $\bar{r}_{obs}(s)$ identifies the radius of the circular boundary at a given arc-length $s$.
For environments with rectangular sections, the inequality \eqref{eq:pos_constr} becomes
\begin{align}
\Big( \frac{2 \bar{w}_i(s) - (\bar{w}_{i,max}(s)+\bar{w}_{i,min}(s))}{(\bar{w}_{i,max}(s)-\bar{w}_{i,min}(s))} \Big)^2 -1 \leq 0,  
\label{eq:rect_constr}
\end{align}
$\forall i = 1,2,$ where $\bar{w}_{i,min}(s)$ and $\bar{w}_{i,max}(s)$ are the
lower and upper bounds at a given arc-length $s$, defining the boundaries of the region.
The constraint boundaries are arc-length functions suitable to model
  fairly complex regions. They represent the physical boundary of a region
  and they can be shaped in order to take into account the presence of
  obstacles attached to the boundary.  As an illustrative example, let us consider the environment with
  rectangular sections depicted in Figure~\ref{fig:w1_obs}. An obstacle
  restricts the collision-free space inside the physical boundary of the
  region.

Let us denote by $\bar{w}^{PB}_{i,max}(s)$ and $\bar{w}^{PB}_{i,min}(s)$ the (respectively) positive and negative distance of the
physical boundary from the frame path at a given arc-length $s$.
We first set $\bar{w}_{i,min}(s) = \bar{w}^{PB}_{i,min}(s)$ and $\bar{w}_{i,max}(s)
= \bar{w}^{PB}_{i,max}(s)$, $ \forall s\in [0,L]$. Then, in order to
take into account the obstacle, we suitably restrict the bounds as follows. 
Let us consider a point $A$ on the
boundary surface of the obstacle.  Let $\boldsymbol{p}_{obs}$ be the position of
point $A$ with components in the inertial frame. We map $\boldsymbol{p}_{obs}$ in
the transverse coordinate vector
$\bar{\boldsymbol{w}}_{obs}$.
First, we select the arc-length on the frame
path, identifying the point at minimum distance from $A$, as 
\begin{align}
s_{obs}:=\text{arg} \min_{s \in \real_0^+}
\|\boldsymbol{p}_{obs}-\bar{\boldsymbol{p}}_f(s)\|^2.
\label{eq:sobs}
\end{align}
Second, we map $\boldsymbol{p}_{obs} - \bar{\boldsymbol{p}}_f(s_{obs})$ into a vector with components in the Serret-Frenet frame attached to the point identified by $s_{obs}$, obtaining
\begin{align}
\bar{w}_{1,obs} &= \bar{\boldsymbol{n}}^T(s_{obs}) (\boldsymbol{p}_{obs} - \bar{\boldsymbol{p}}_f(s_{obs})), \label{eq:w1obs}\\
\bar{w}_{2,obs} &= \bar{\boldsymbol{b}}^T(s_{obs}) (\boldsymbol{p}_{obs} - \bar{\boldsymbol{p}}_f(s_{obs})). \label{eq:w2obs}
\end{align}

Since, according to the particular scenario, the obstacle only affects the function $\bar{w}_{1,min}(\cdot)$, we update
\begin{align*}
\bar{w}_{1,min}(s_{obs}) &= \max\{\bar{w}^{PB}_{1,min}(s_{obs}),
\bar{w}_{1,obs}\}.
\end{align*}

\subsection{Equivalent minimum-time formulation and optimal control solver}
\label{sec:mintime-new}
The minimum-time problem \eqref{eq:mintime_standard} is reformulated in the new (arc-length dependent) variables $\bar{\boldsymbol{x}}_w$ and $\bar{\boldsymbol{u}}$, by using the cost \eqref{eq:cost_s}, the transverse dynamics \eqref{eq:tran_dynamics} and the constraints \eqref{eq:cT}, \eqref{eq:constr_s}, and \eqref{eq:pos_constr}. Denoting by $c(\bar{\boldsymbol{x}}_w(s),\bar{\boldsymbol{u}}(s)) \leq 0, \; \forall s \in [0,L],$ the constraints \eqref{eq:constr_s} and \eqref{eq:pos_constr} in vectorial form, the reformulated problem is
\begin{equation}
  \begin{split}
   \label{eq:mintime}
    \min_{\bar{\boldsymbol{x}}_w(\cdot),\bar{\boldsymbol{u}}(\cdot)} &\;   \int_0^L \!\!\! \quad \frac{1-\bar{k}(s) \bar{w}_1(s)}
 {\bar{\boldsymbol{t}}(s)^T \bar{\text{\textbf{v}}}(s)}
 \; ds,\\
    \!\!\subj &\;
    \bar{\boldsymbol{x}}'_w (s)= \bar{f}(\bar{\boldsymbol{x}}_w(s),\bar{\boldsymbol{u}}(s)), \quad \bar{\boldsymbol{x}}_w(0) = \boldsymbol{x}_{w0},\\[0.5ex]
    & \; c_f(\bar{\boldsymbol{x}}_w(L)) \leq 0,\\
    & c(\bar{\boldsymbol{x}}_w(s),\bar{\boldsymbol{u}}(s))\leq 0, \; \forall s \in [0,L].
  \end{split}
\end{equation}
Note that the fixed horizon problem \eqref{eq:mintime} is equivalent to \eqref{eq:mintime_standard} since
trajectories solving \eqref{eq:mintime} can be mapped into trajectories solving
\eqref{eq:mintime_standard}. 

In order to solve problem \eqref{eq:mintime}, we use a combination of the
PRojection Operator based Newton method for Trajectory Optimization (PRONTO)
\cite{hauser2002projection} with a barrier function approach
\cite{hauser2006barrier}. 
We relax state-input constraints by adding them in
the cost functional, i.e., we consider the problem
\begin{equation}
\begin{split}
\min_{\bar{\boldsymbol{x}}_w(\cdot),\bar{\boldsymbol{u}}(\cdot)} \!\!\!\!&\;\; \int_0^L \!\!\! 
\Big( 
\frac{1-\bar{k}(s) \bar{w}_1(s)}{\bar{\boldsymbol{t}}(s)^T \bar{\text{\textbf{v}}}(s)} + \epsilon \sum_j \beta_\nu (-c_j(\bar{\boldsymbol{x}}_w(s),\bar{\boldsymbol{u}}(s))) 
\Big) ds\\
& \quad + \; \epsilon_f \sum_i \beta_{\nu_f} (-c_{f,i}(\bar{\boldsymbol{x}}_w(L))),\\
\!\!\subj &\;\; \bar{\boldsymbol{x}}'_w(s) = \bar{f}(\bar{\boldsymbol{x}}_w(s),\bar{\boldsymbol{u}}(s)), \quad \forall s \in [0,L],\\
& \;\; \bar{\boldsymbol{x}}_w(0) = \boldsymbol{x}_{w0}.
\end{split}
\label{eq:mintime2}
\end{equation}
where $\epsilon$ and $\epsilon_f$ are positive parameters and $\beta_\ell(\cdot)$,
  $\ell \in \{\nu,\nu_f\}$, is a function depending on the parameter $\ell$ and
  defined as
\begin{align*}                
\beta_\ell(x) &:=  
                  \begin{cases} 
                    -\log(x) & x > \ell,\\
                    -\log(\ell) + \frac{1}{2} \big[ (\frac{x - 2\ell}{\ell})^2 -1
                    \big] & x \leq \ell.
                  \end{cases} \nonumber           
\end{align*}

Let an initial trajectory for the initialization of the algorithm be given.
The strategy to find an approximated solution to \eqref{eq:mintime} can be summarized as follows.
Problem \eqref{eq:mintime2} is iteratively solved
by reducing the parameters $\epsilon, \nu, \epsilon_f$ and $\nu_f$ at each iteration, and thus pushing
the trajectory towards the constraint boundaries.
Each instance of problem \eqref{eq:mintime2} is solved by means of the PRONTO
algorithm described in Appendix \ref{app:pronto}.

\subsection{Summary of the strategy}
\label{sec:summary}
A pseudo code of the whole strategy to compute minimum time trajectories 
is reported in the following table (Algorithm 1). We denote with 
$(\bar{\boldsymbol{x}}_{w}(\cdot),\bar{\boldsymbol{u}}(\cdot))^0$ the initial trajectory to initialize the algorithm and with
$\ponewt$ the PRojection Operator based Newton method for Trajectory Optimization routine
that, given a trajectory $(\bar{\boldsymbol{x}}_{w}(\cdot),\bar{\boldsymbol{u}}(\cdot))^{i-1}$, computes the solution $(\bar{\boldsymbol{x}}_{w}(\cdot),\bar{\boldsymbol{u}}(\cdot))^i$ to problem \eqref{eq:mintime2}, i.e., $(\bar{\boldsymbol{x}}_{w}(\cdot),\bar{\boldsymbol{u}}(\cdot))^i = \ponewt((\bar{\boldsymbol{x}}_{w}(\cdot),\bar{\boldsymbol{u}}(\cdot))^{i-1})$.
\begin{algorithm}[H]
\caption{Minimum-time strategy}
\label{alg:proj_newt}
\begin{algorithmic}
\REQUIRE initial condition $\boldsymbol{x}_0$, final desired region $X_T$, bounds $p_{max}, q_{max}, r_{max}, f_{min}, f_{max}, \varphi_{max}(\!\cdot\!),
\theta_{max}\!(\!\cdot\!),$ $\psi_{max}\!(\!\cdot\!)$, and the dynamic model \eqref{eq:state_space}
\STATE
\textbf{A. Frame path} \\
generate $\bar{\boldsymbol{p}}_f(s), \; \forall s \in [0,L]$\\
compute 
\begin{itemize} 
\item tangent, normal and binormal vectors\\[0.5ex]
$\bar{\boldsymbol{t}}(s)\!=\!\bar{\boldsymbol{p}}'_f(s)$, \;
$\bar{\boldsymbol{n}}(s)\!=\!\frac{\bar{\boldsymbol{p}}''_f(s)}{\Vert \bar{\boldsymbol{p}}''_f(s)\Vert_2}$, \;
$\bar{\boldsymbol{b}}(s)\!=\!\bar{\boldsymbol{t}}(s) \times \bar{\boldsymbol{n}}(s)$
\item curvature
$\bar{k}(s)\!=\!\Vert \bar{\boldsymbol{p}}''_f(s)\Vert_2$
\item torsion
$\bar{\tau}(s)\!=\!\bar{\boldsymbol{n}}(s) \; \bar{\boldsymbol{b}}'(s)$
\end{itemize}
\textbf{B. Transverse dynamics} \\
\STATE set-up transverse dynamics \eqref{eq:transv-all} \\[0.5ex]
\textbf{C. Cost and constraints} \\[0.5ex]
\STATE
set-up cost $ \int_0^L \frac{1-\bar{k}(s) \bar{w}_1(s)}{\bar{\boldsymbol{t}}(s)^T \bar{\text{\textbf{v}}}(s)} \; ds$\\[0.5ex]
set-up constraints \eqref{eq:cT} and \eqref{eq:constr_s}\\[0.5ex]
define boundary constraints by using \eqref{eq:circ_constr} and/or \eqref{eq:rect_constr} \\[0.5ex]

\textbf{E. Numerical solution to \eqref{eq:mintime}} \\
\STATE
compute initial trajectory $(\bar{\boldsymbol{x}}_{w}(\cdot),\bar{\boldsymbol{u}}(\cdot))^0$\\
\STATE 
set $\epsilon = 1$, $\epsilon_f = 1$, $\nu = 1$, $\nu_f = 1$     
\FOR{$i = 1, 2 \ldots$}
\STATE compute: $\!(\bar{\boldsymbol{x}}_{w}(\cdot),\bar{\boldsymbol{u}}(\cdot))^i \!=\! \ponewt((\bar{\boldsymbol{x}}_{w}(\cdot),\bar{\boldsymbol{u}}(\cdot))^{i-1})$
\STATE update: $\epsilon$, $\epsilon_f$, $\nu$, $\nu_f$
\ENDFOR
\ENSURE $(\bar{\boldsymbol{x}}_{w}(\cdot),\bar{\boldsymbol{u}}(\cdot))_{opt} =(\bar{\boldsymbol{x}}_{w}(\cdot),\bar{\boldsymbol{u}}(\cdot))^{i}$
\end{algorithmic}
\end{algorithm}

\section{Numerical computations}
\label{sec:computations}
In this section, we present numerical computations and experimental tests on a nano-quadrotor with mass $m = 0.0325 \; \text{kg}$, in order to show the effectiveness of the proposed strategy.
First, we consider a scenario with two obstacles: a parallelepiped and a cylinder, as depicted in Figure \ref{fig:path_s2}.
Second, we consider an experimental scenario 
and we show the results related to the execution of the optimal trajectory using our maneuver regulation control scheme \cite{SS-GN-HHB-AF:13}.
\subsection{Rooms with obstacles}
\begin{figure}[htbp]
		\begin{center}
			\vspace{-0.4cm}\subfloat[][Path: 3D view]
			{\includegraphics[scale=0.47]{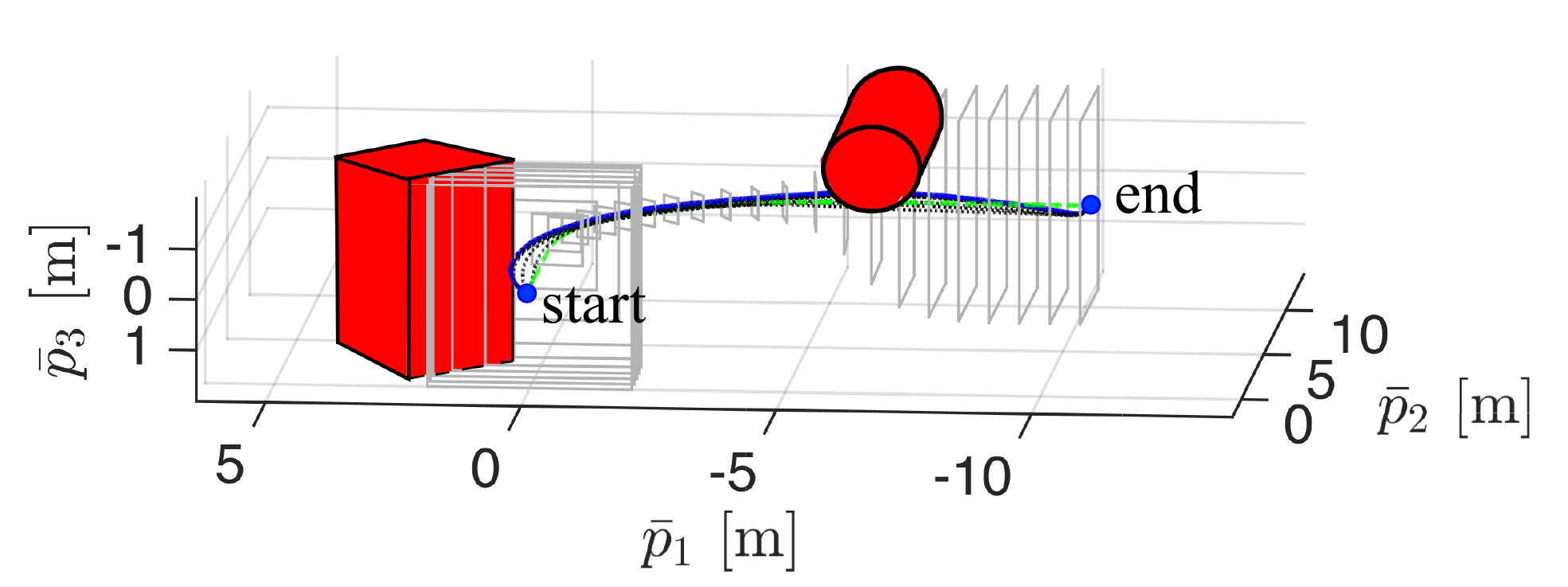}\label{fig:path_s2}}\\
			\vspace{-0.2cm}
			\subfloat[][Transverse coordinate $\bar{w}_1$]
			{\hspace{-0.1cm}\includegraphics[width=4.7cm]{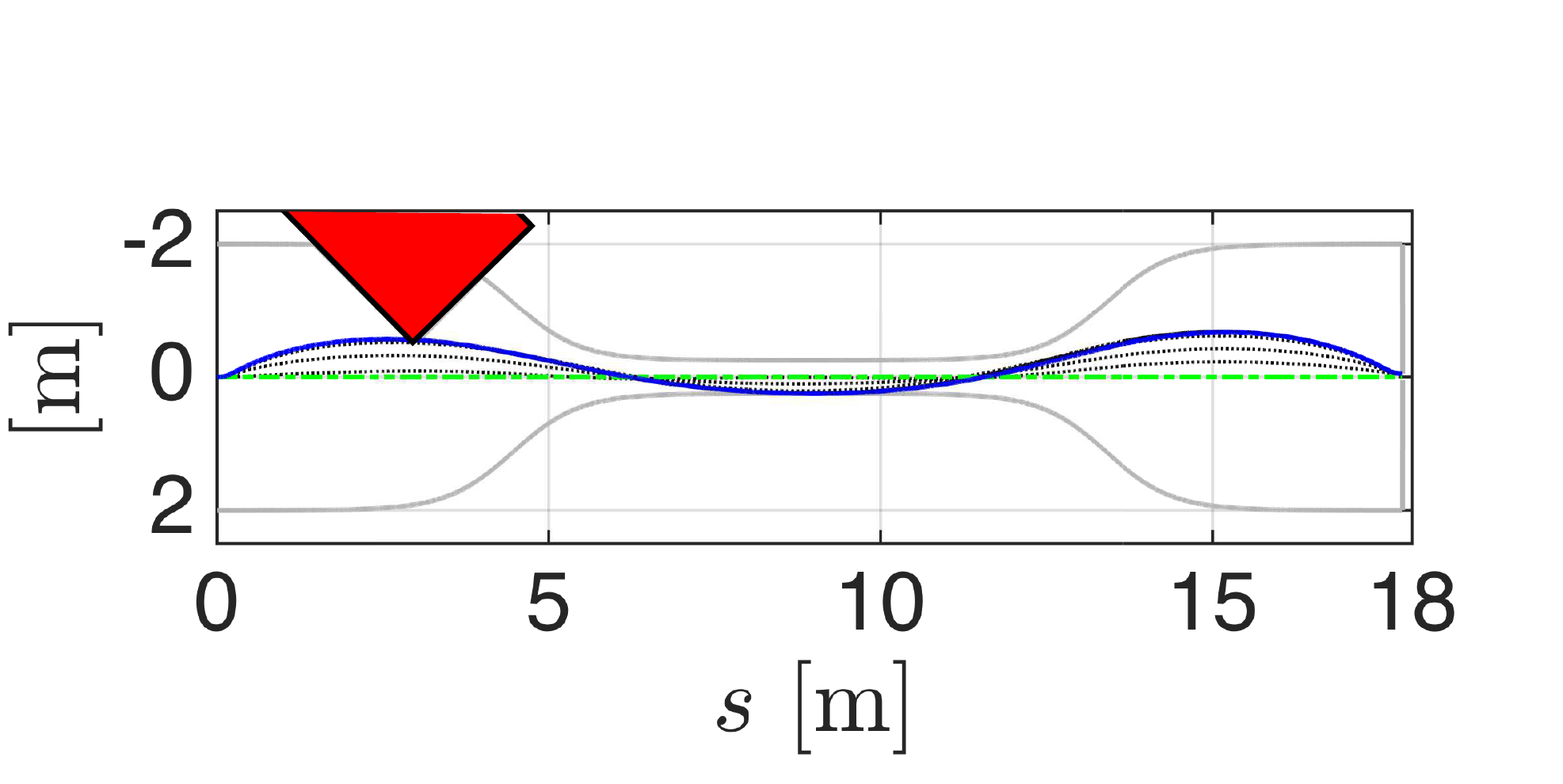}\label{fig:w1_s2}}
			\subfloat[][Transverse coordinate $\bar{w}_2$]
			{\vspace{-0.5cm}\includegraphics[width=4.7cm]{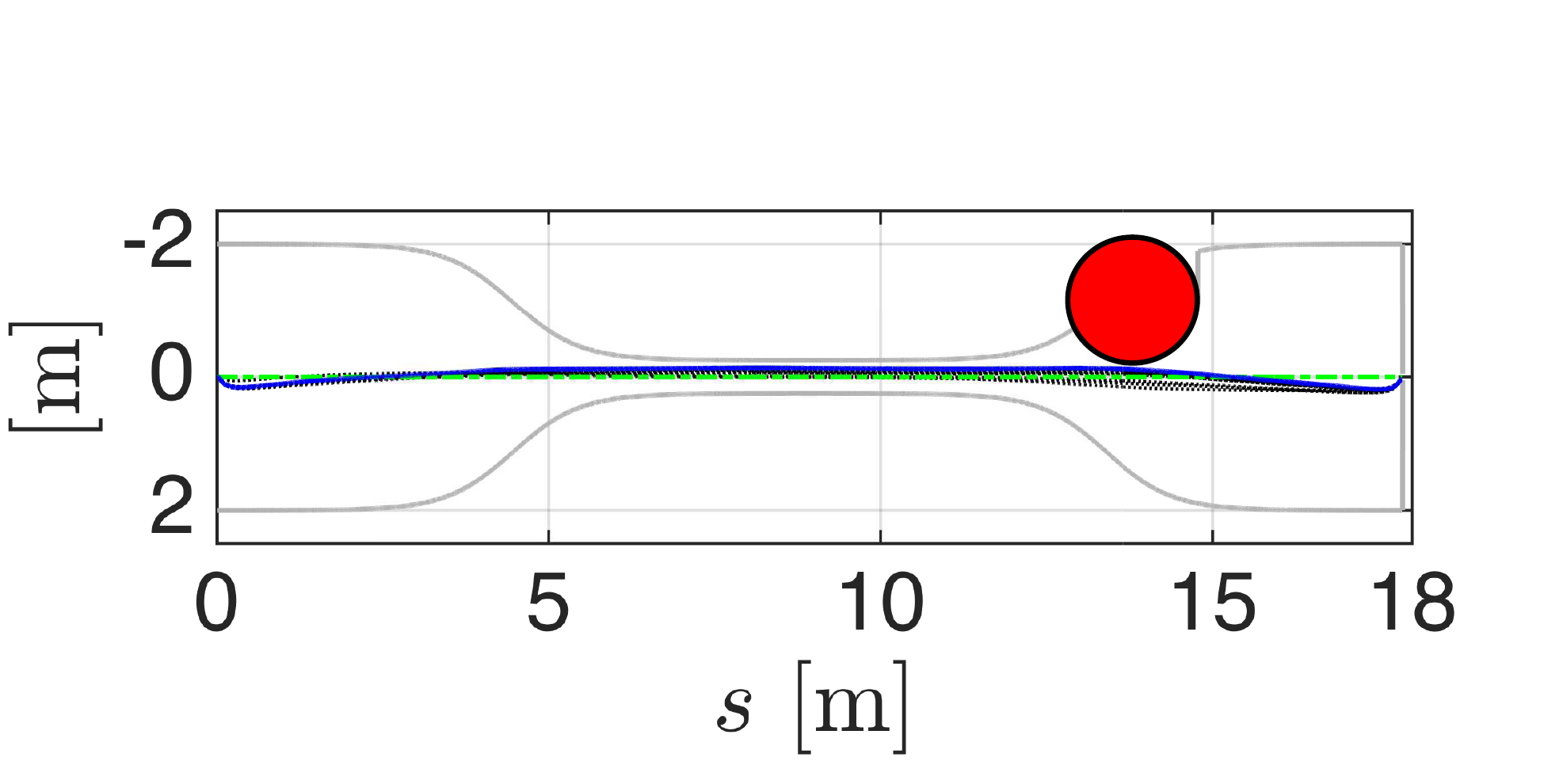}\label{fig:w2_s2}}\\
			\caption{Path and transverse coordinates. Initial (dot-dashed green), intermediate (dotted black) and minimum-time (solid blue) trajectory. Constraint boundaries are depicted in grey.}
			\label{fig:scenario2_path}
		\end{center}
\end{figure}

\begin{figure}[htbp]
			\begin{center}
			\vspace{-0.5cm}\subfloat[][Velocity $\bar{\boldsymbol{t}}^T \bar{\text{\textbf{v}}}$]
			{\includegraphics[width=4.7cm]{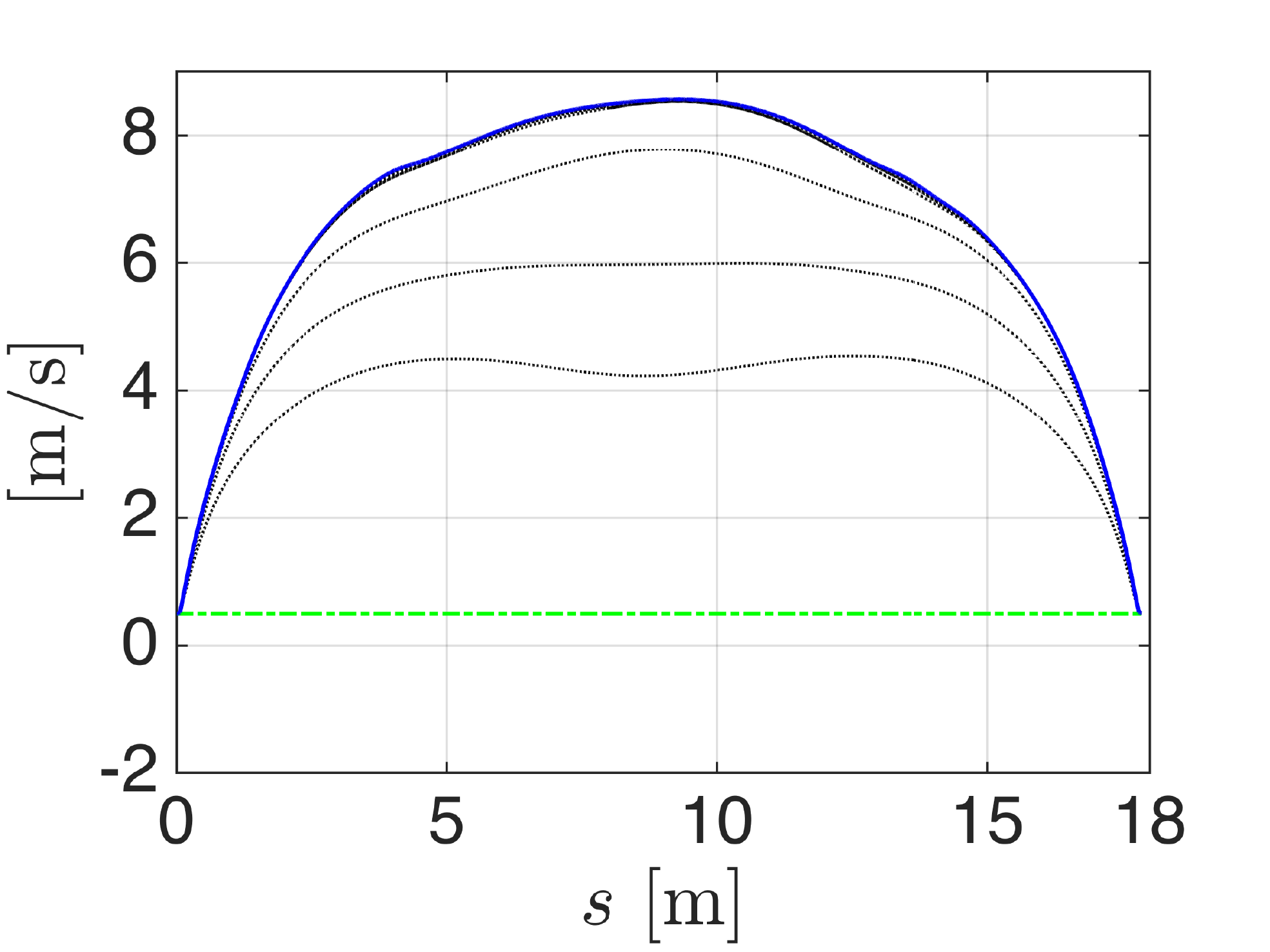}\label{fig:dp1_s2}}
			\subfloat[][Roll angle $\bar{\varphi}$]
			{\hspace{-0.25cm}\includegraphics[width=4.7cm]{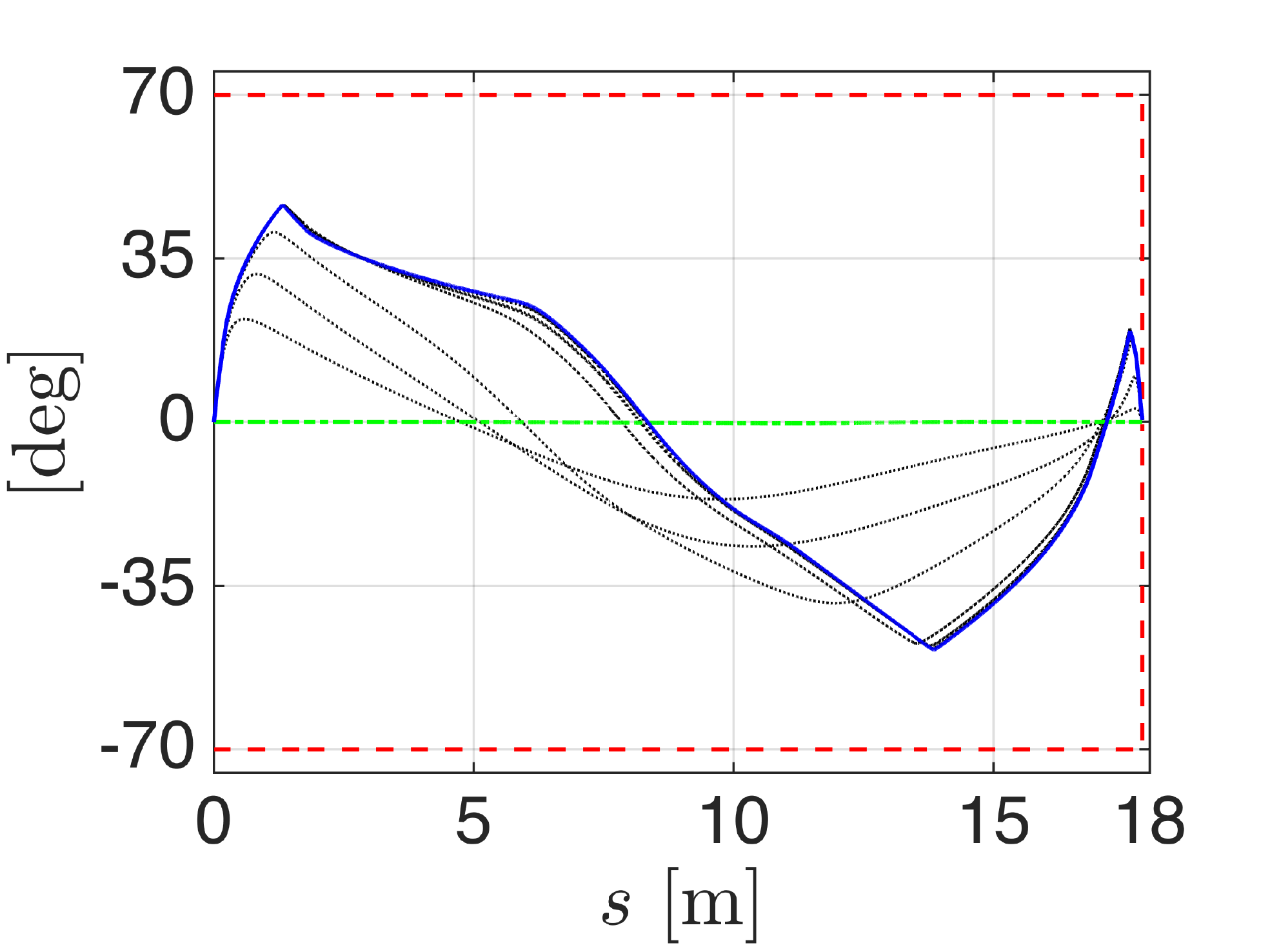}\label{fig:phi_s2}}\\
			\subfloat[][Velocity $\bar{\boldsymbol{n}}^T \bar{\text{\textbf{v}}}$]
			{\includegraphics[width=4.7cm]{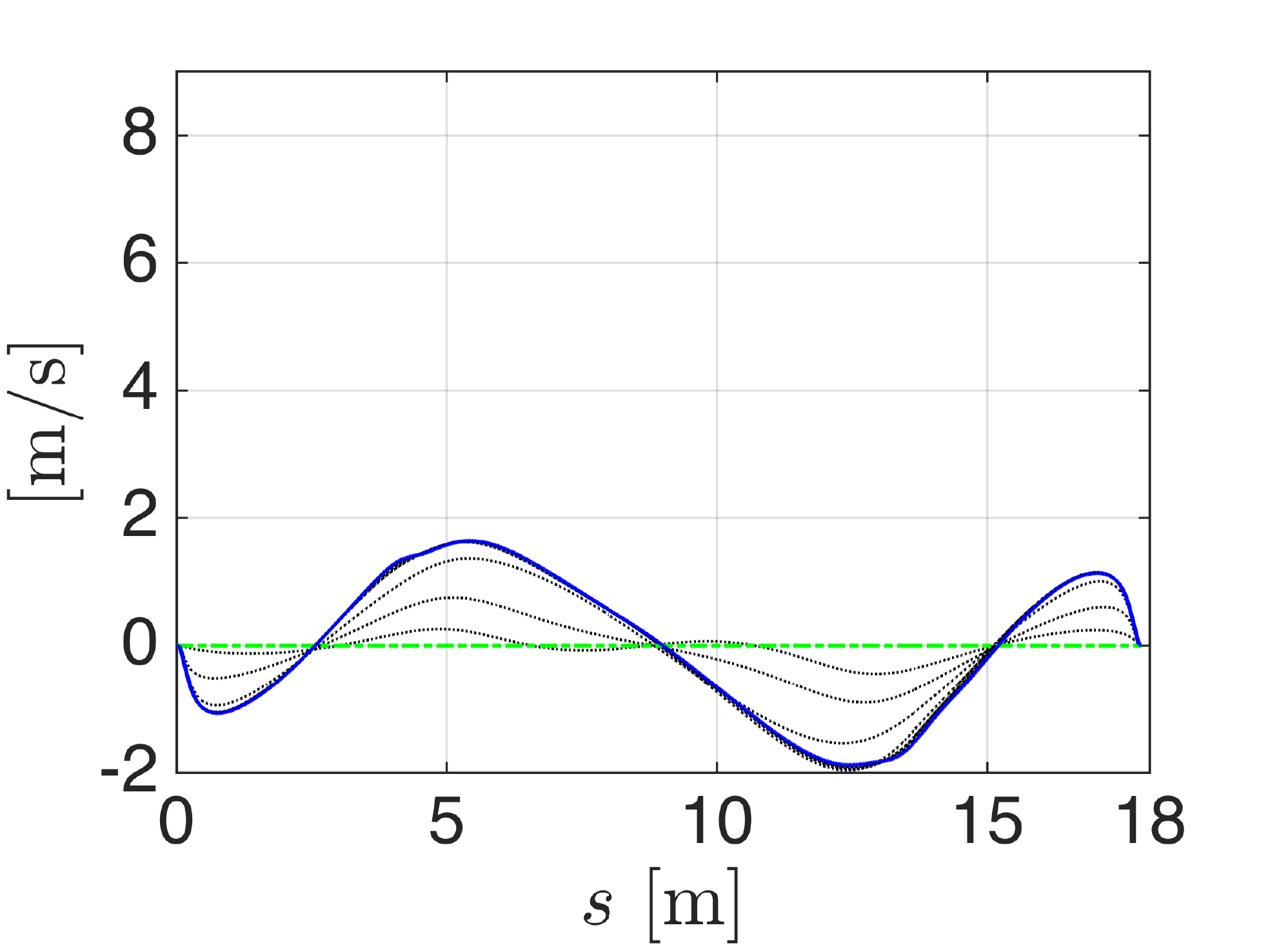}\label{fig:dp2_s2}}
			\subfloat[][Pitch angle $\bar{\theta}$]
			{\hspace{-0.25cm}\includegraphics[width=4.7cm]{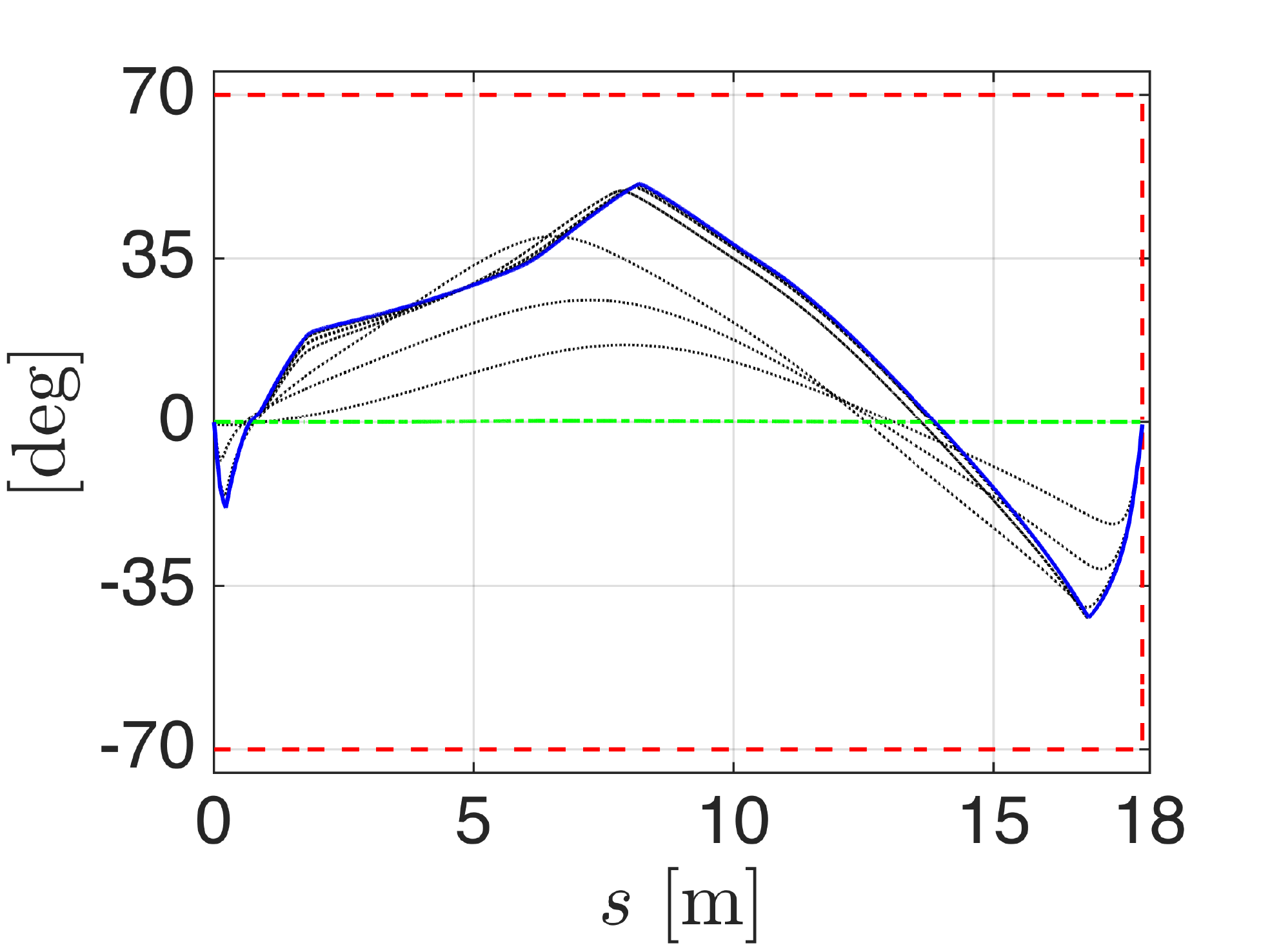}\label{fig:th_s2}}\\
			\subfloat[][Velocity $\bar{\boldsymbol{b}}^T \bar{\text{\textbf{v}}}$]
			{\includegraphics[width=4.7cm]{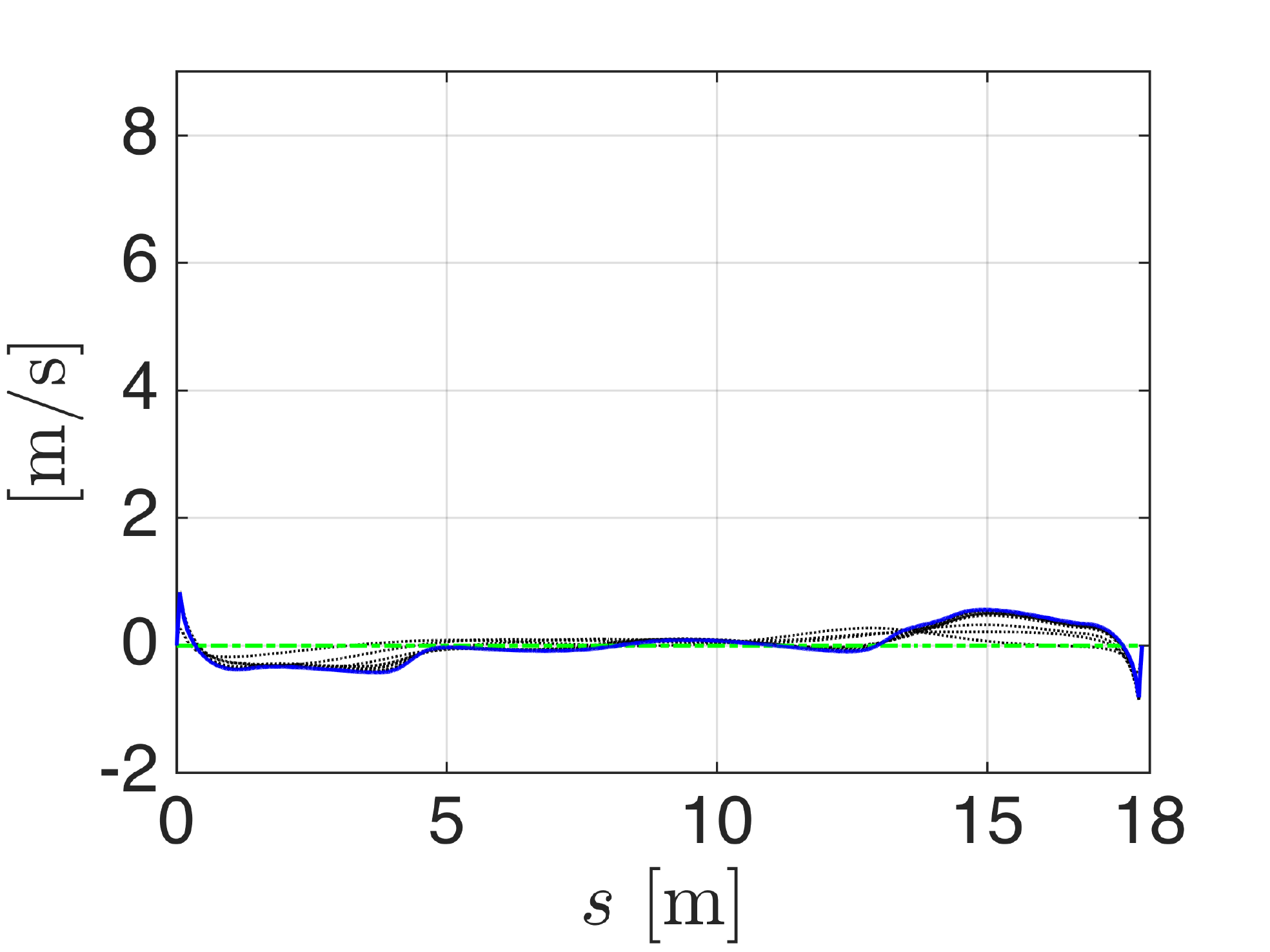}\label{fig:dp3_s2}}
			\subfloat[][Yaw angle $\bar{\psi}$]
			{\hspace{-0.25cm}\includegraphics[width=4.7cm]{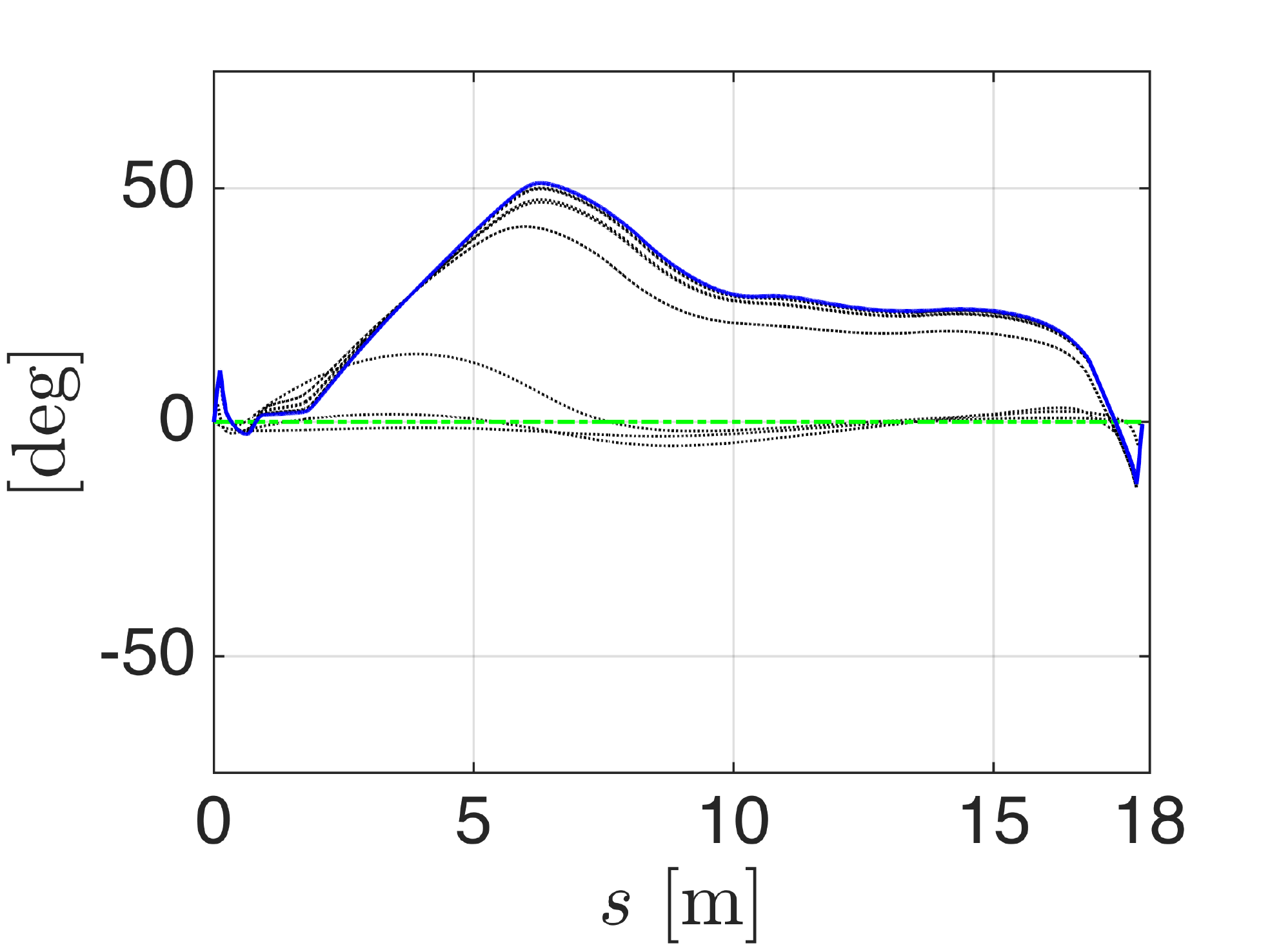}\label{fig:psi_s2}}\\							
			\caption{Velocities and angles. Initial (dot-dashed green), intermediate (dotted black) and minimum-time (solid blue) trajectory. Constraint boundaries are depicted in red.}
			\label{fig:scenario2_velang}
		\end{center}
	\end{figure}

\begin{figure}[htbp]
		\begin{center}
			\vspace{-0.7cm}\subfloat[][Roll rate $\bar{p}$]
			{\includegraphics[width=4.7cm]{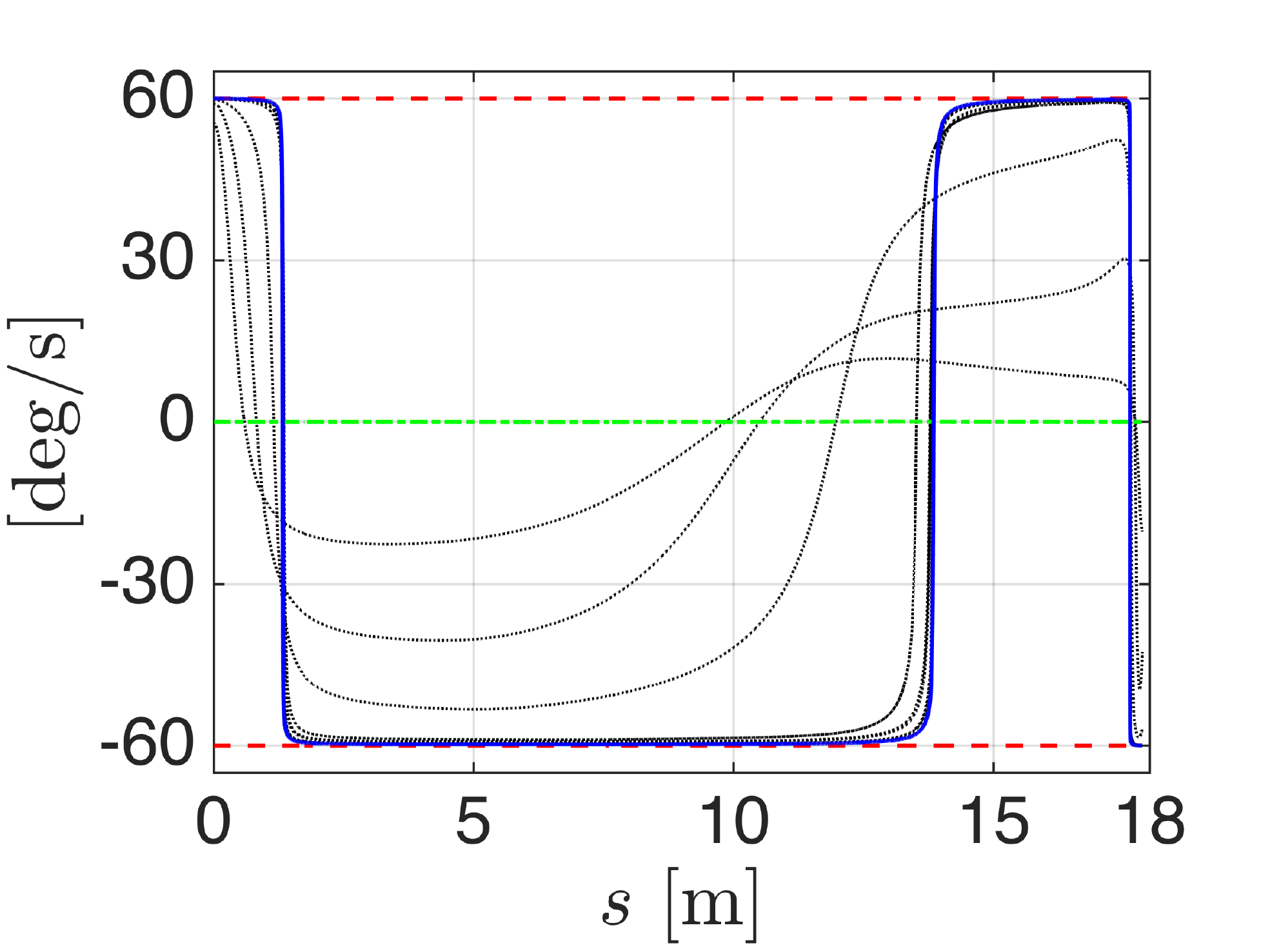}\label{fig:pp_s2}}
			\subfloat[][Pitch rate $\bar{q}$]
			{\hspace{-0.25cm}\includegraphics[width=4.7cm]{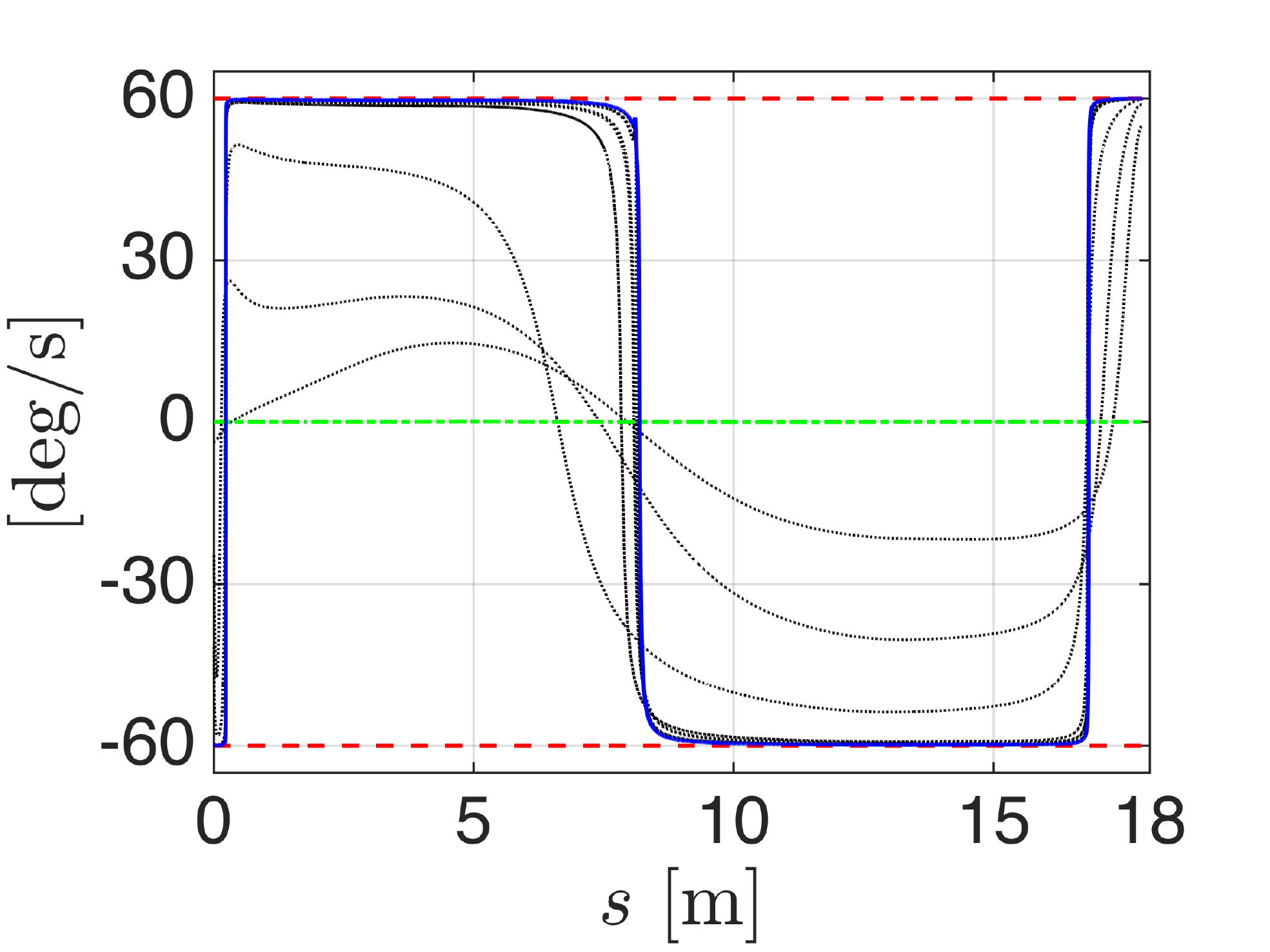}\label{fig:qq_s2}}\\
			\subfloat[][Yaw rate $\bar{r}$]
			{\includegraphics[width=4.7cm]{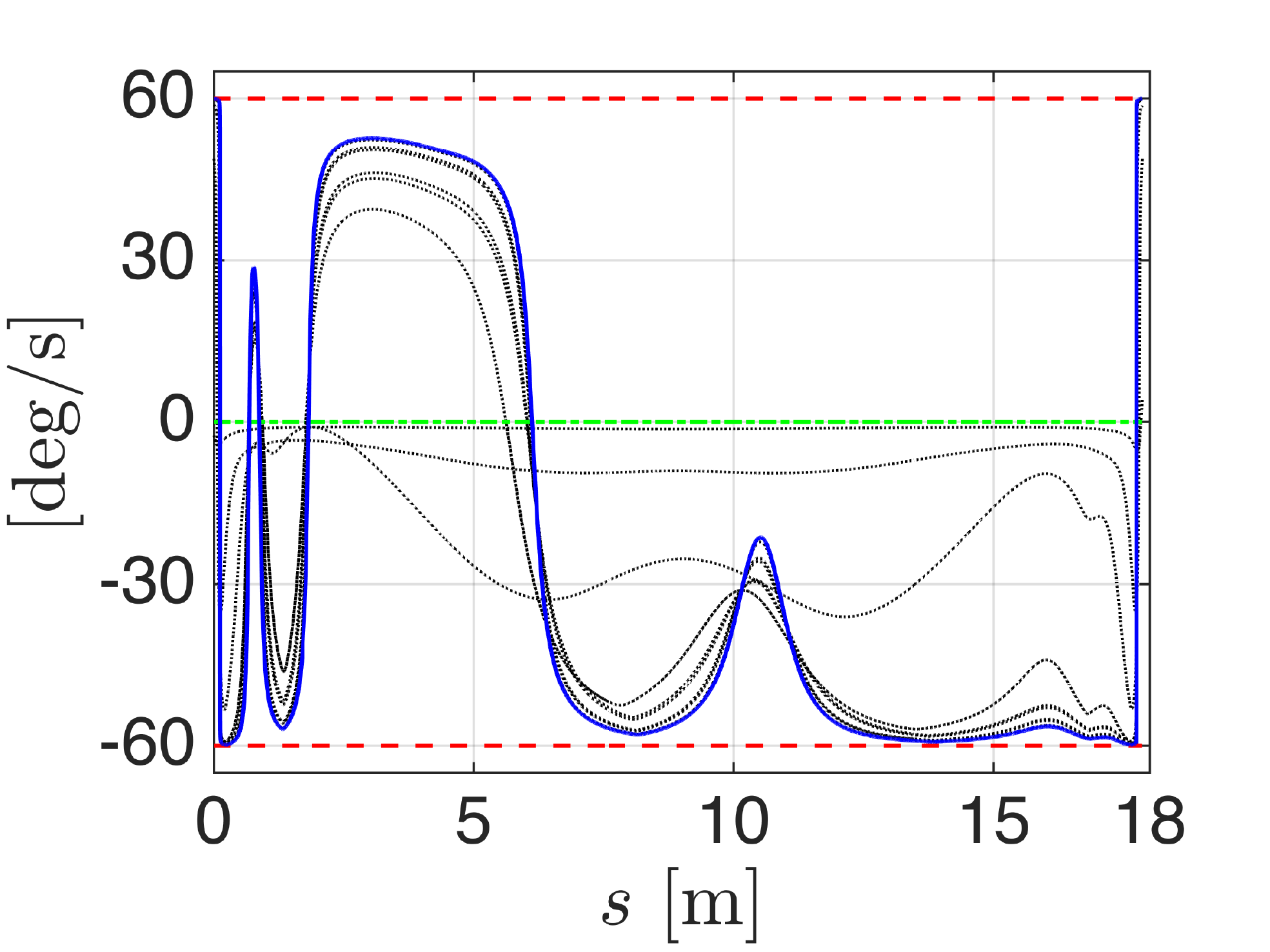}\label{fig:rr_s2}}
			\subfloat[][Thrust $\bar{f}$]
			{\hspace{-0.25cm}\includegraphics[width=4.7cm]{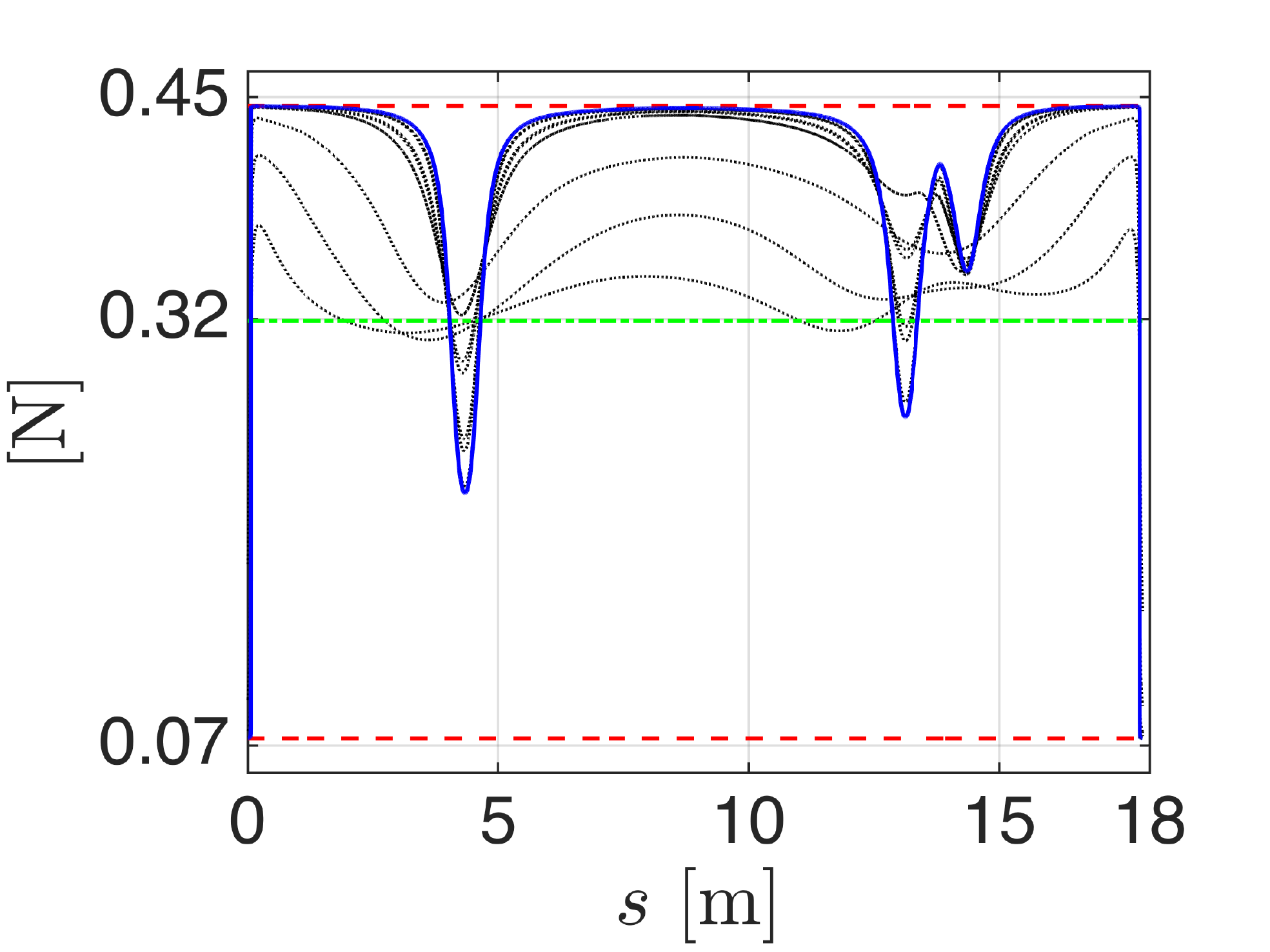}\label{fig:ff_s2}}\\						
			\caption{Inputs. Initial (dot-dashed green), intermediate (dotted black) and minimum-time (solid blue) trajectory. Constraint boundaries are depicted in dashed red.}
			\label{fig:scenario2_inputs}
		\end{center}
	\end{figure}

The first scenario is as follows. The vehicle has to move from one room to another through a
narrow corridor.
There is a parallelepiped in the first room and a cylinder in the second room.
As an additional
requirement, the quadrotor must reach a neighborhood of $\boldsymbol{x}_{w0}$ at the end of its motion.
In order to fulfil this objective, we consider the final constraint \eqref{eq:cT} with $  c_{f,i}(\bar{x}_{w_i}(L))$ as in \eqref{eq:c_fi}, where
$\bar{x}_{w_i, min} = x_{w_i,0} - \text{tol}_i$, $\bar{x}_{w_i, max} = x_{w_i,0} + \text{tol}_i$, $\text{tol}_i$ is a given tolerance and $x_{w_i,0}$ is the $i$-th component of $\boldsymbol{x}_{w0}$.
Results are depicted in
Figures \ref{fig:scenario2_path}, \ref{fig:scenario2_velang},
\ref{fig:scenario2_inputs}.  The initial trajectory is depicted in dot-dashed
green, intermediate trajectories in dotted black and the minimum-time trajectory
in solid blue. Collision-free boundaries are depicted in grey and remaining
state-input constraints in dashed red.

We choose as frame path a $C^{\infty}$ curve on the $\bar{p}_1 - \bar{p}_2$
plane with constant binormal vector $\bar{\boldsymbol{b}} = [0 \; 0 \; 1]^T$ and curvature
\begin{align*}
\bar{k}(s) = \frac{1}{5} \; \frac{\tanh(s-5) - \tanh(s-5(1+\frac{\pi}{2}))}{\max(\tanh(s-5) - \tanh(s-5(1+\frac{\pi}{2})))}.
\end{align*}

The collision free region is defined by constraint \eqref{eq:rect_constr} where
obstacle boundaries $\bar{w}_{i,min}(\cdot)$ and $\bar{w}_{i,max}(\cdot)$,
$i = 1,2$, are chosen as follows. 
Functions $\bar{w}_{1,max}$ and $\bar{w}_{2,max}$ are not affected by obstacles.
As depicted in Figures \ref{fig:w1_s2} and \ref{fig:w2_s2},
$\bar{w}_{1,max}(\cdot)$ and $\bar{w}_{2,max}(\cdot)$ are obtained using sigmoid functions
with values varying from $2 \; \text{m}$ to $0.25 \; \text{m}$.
Functions  $\bar{w}_{1,min}$ and $\bar{w}_{2,min}$ are affected by obstacles.
In order to model the obstacles, we consider the position of the obstacle boundary as a function of its arc-length.
We choose $10^{-3}$ as discretization step for the arc-length 
and for every value of the boundary position $\boldsymbol{p}_{obs}$, we compute $s_{obs}$ and $\bar{w}_{1,obs}, \bar{w}_{2,obs}$ by using equations \eqref{eq:sobs}, \eqref{eq:w1obs} and \eqref{eq:w2obs}, respectively.
Thus, in order to define $\bar{w}_{1,min}$ and $\bar{w}_{2,min}$, we
first set $\bar{w}_{1,min}(s) = -\bar{w}_{1,max}(s)$ and $\bar{w}_{2,min}(s) = -\bar{w}_{2,max}(s)$, 
$\forall s \in [0,L]$. 
Second, for each $s_{obs}^R$ and $\bar{w}_{1,obs}^R$ related to a point $R$ on the parallelepiped, we update 
\begin{align*}
\bar{w}_{1,min}(s_{obs}^R) &= \max\{-\bar{w}_{1,max}(s_{obs}^R), \bar{w}_{1,obs}^R\}.
\end{align*}
Third and final, for each $s_{obs}^C$ and $\bar{w}_{2,obs}^C$ related to a point $C$  on the cylinder,
we update  
\begin{align*}
\bar{w}_{2,min}(s_{obs}^C) &= \max\{-\bar{w}_{2,max}(s_{obs}^C), \bar{w}_{2,obs}^C\}.
\end{align*}

We choose the initial trajectory as follows. We set the frame path as the position, a velocity module of $0.5 \; \text{m/s}$ along the curve and a zero yaw angle.  The remaining initial states and inputs are
computed by using the differential flatness of the quadrotor dynamics \cite{mellinger2011minimum}. It is worth noting that the position part of the initial trajectory does not have to necessarily match the frame path. Also, the initial trajectory could be alternatively computed through the projection of a state-input curve by using the projection operator \eqref{eq:proj_oper_def} described in Appendix \ref{app:pronto}, instead of using the differential flatness.

Having the initial trajectory in hand, we run the algorithm to numerically
  compute solutions. Note that the PRONTO method (described in Appendix
  \ref{app:pronto}) is designed considering an $s$-dependent continuous
  dynamics. In order to implement it by using a numerical toolbox (Matlab), we
  consider a suitable tolerance. We choose $10^{-3}$ as discretization step on
  $s$ and we use the tolerance of the Matlab solver to integrate the
  differential equations. Each intermediate optimal trajectory is computed by
solving the optimization problem \eqref{eq:mintime2} with constant values
of the parameters $\epsilon$, $\nu$, $\epsilon_f$ and $\nu_f$.  We start
with $\epsilon = 1$, $\nu = 1$, $\epsilon_f = 1$, $\nu_f = 1$ and,
  following a suitable heuristic, we decrease them at each iteration.
Since the algorithm
operates in an interior point fashion, intermediate trajectories are all
feasible and are pushed to the constraint boundaries when $\epsilon,\nu,\epsilon_f, \nu_f$ are decreased.

As regards the minimum-time
trajectory, the maneuver is performed in $3.57 \; \text{s}$ and the path touches
the constraint boundaries when the vehicle is inside the corridor and in the proximity of obstacles (Figures
\ref{fig:w1_s2} and \ref{fig:w2_s2}). The velocity $\bar{\boldsymbol{t}}^T \bar{\text{\textbf{v}}}$ (Figure
\ref{fig:dp1_s2}) reaches a peak of about $8.5 \; \text{m/s}$ in the middle of the path
and approaches the final desired value at the end. Velocities $\bar{\boldsymbol{n}}^T \bar{\text{\textbf{v}}}$ and
$\bar{\boldsymbol{b}}^T \bar{\text{\textbf{v}}}$ (Figures \ref{fig:dp2_s2} and \ref{fig:dp3_s2}, respectively)
are between $-2.0 \; \text{m/s}$ and $2.0 \; \text{m/s}$. Roll and pitch angles (Figures \ref{fig:phi_s2}, \ref{fig:th_s2},
respectively) do not touch constraint boundaries and alternate positive and
negative values between $-50 \; \text{deg}$ and $50 \; \text{deg}$. 
The yaw angle has values between $-20 \; \text{deg}$ and $50 \; \text{deg}$ (Figure \ref{fig:psi_s2}).
As regards
the inputs, while constraints on roll and pitch rates (Figures \ref{fig:pp_s2},
\ref{fig:qq_s2}, respectively) are always active, yaw rate and thrust (Figures
\ref{fig:rr_s2} and \ref{fig:ff_s2}, respectively) alternate intervals with
active and inactive constraints. Furthermore, note that the final state reaches a neighborhood 
of the initial state, satisfying $|| \bar{\boldsymbol{x}}_w(L) - \boldsymbol{x}_{w0}|| < 0.07$.

\subsection{{Tubular passage}}	
As a second test, 
we consider a region delimited by hula hoops as constrained environment.
First, we compute a minimum-time trajectory through our optimization strategy and second, we 
experimentally execute the minimum-time trajectory on the CrazyFlie nano-quadrotor (https://www.bitcraze.io/crazyflie/), by using a suitable controller.
We invite the reader to watch the attached video related to this experiment.

We set-up the optimization algorithm as follows.
We approximate the collision free region as a tube with circular section.
We choose as {frame path} a curve on the $\bar{p}_2 - \bar{p}_3$
plane with constant binormal vector $\bar{\boldsymbol{b}} = [1 \; 0 \; 0]^T$ and curvature
\begin{align*}
\bar{k}(s) = \frac{\frac{1}{1+e^{-8(s-2.27)}}-\frac{1}{1+e^{-8(s-3.67)}}}{\max(\frac{1}{1+e^{-8(s-2.27)}}-\frac{1}{1+e^{-8(s-3.67)}})}.
\end{align*}
Moreover, we consider the constraint {\eqref{eq:circ_constr}} with {constant} $\bar{r}_{obs}=r_{hh}-l-e_p$, where $r_{hh} = 0.33 \; \text{m}$ is the hula hoop radius, $l = 0.04 \; \text{m}$ is the distance between the quadrotor center of mass and its propellers and $e_p = 0.01 \; \text{m}$ is the estimated position error arising during control. 

As regards input constraints, we impose, for safety reasons, more severe bounds than the ones required by the physical vehicle limitations. In this way, we also ensure that the ``experimental" trajectory remains feasible although the imperfect tracking of desired inputs by actual values (naturally arising during control).
We choose $p_{min}=-15 \; \text{deg/s}$ and $p_{max}=15 \; \text{deg/s}$ for the roll rate and $\force_{min}= 0.1779 \; \text{N}$ and $\force_{max}= 0.3411 \; \text{N}$ for thrust. 

\begin{figure}[htbp]
\hspace{-0.5cm}
	\begin{tabular}{cc}
		{\multirow{-2}[15.0]{*}{\subfloat[][Path]{\includegraphics[width=4.3cm]{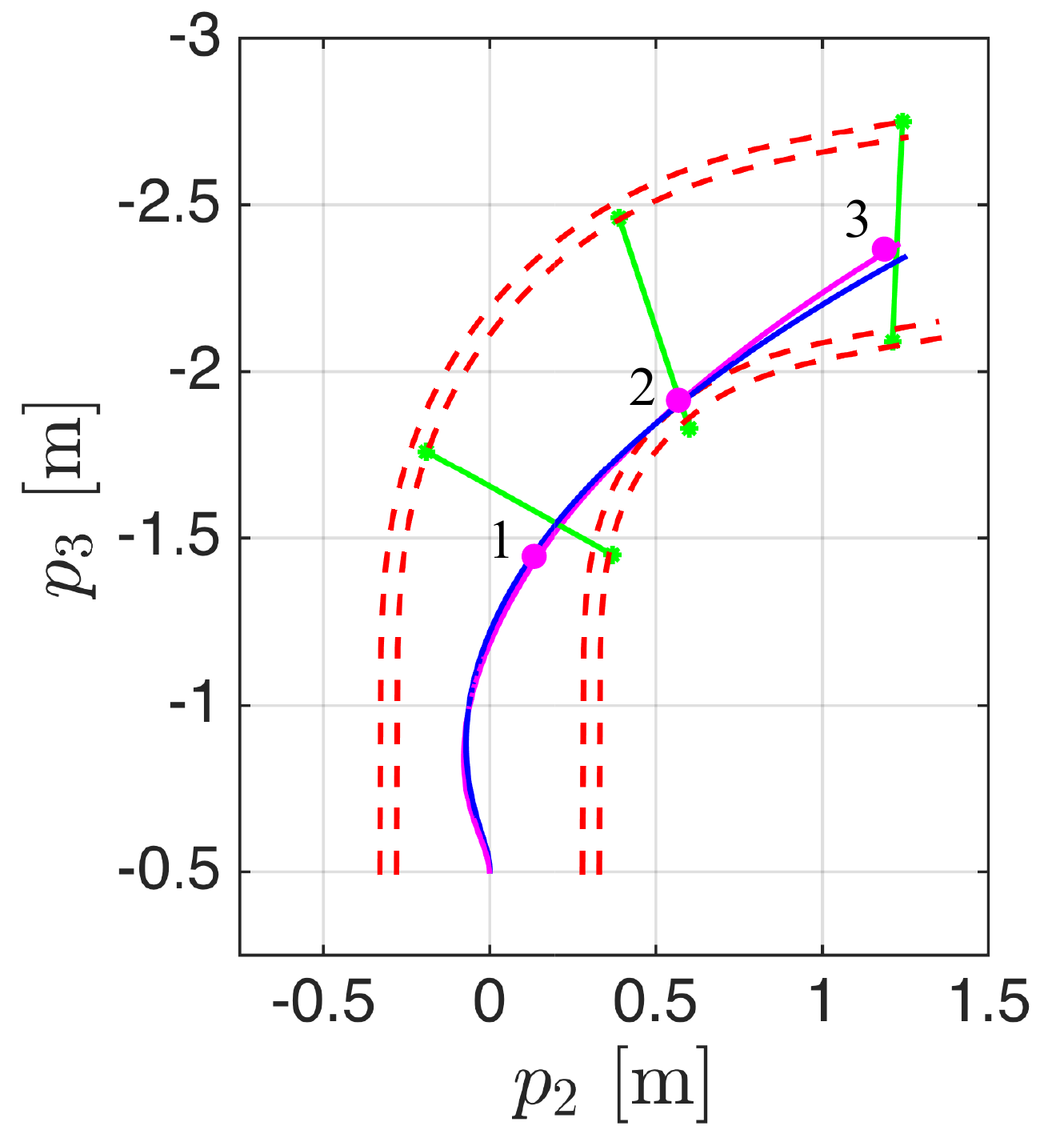} \label{fig:path_exp}}}} &
		\hspace{-0.7cm}{\subfloat[][Roll angle $\varphi$]{{\includegraphics[width=4.7cm]{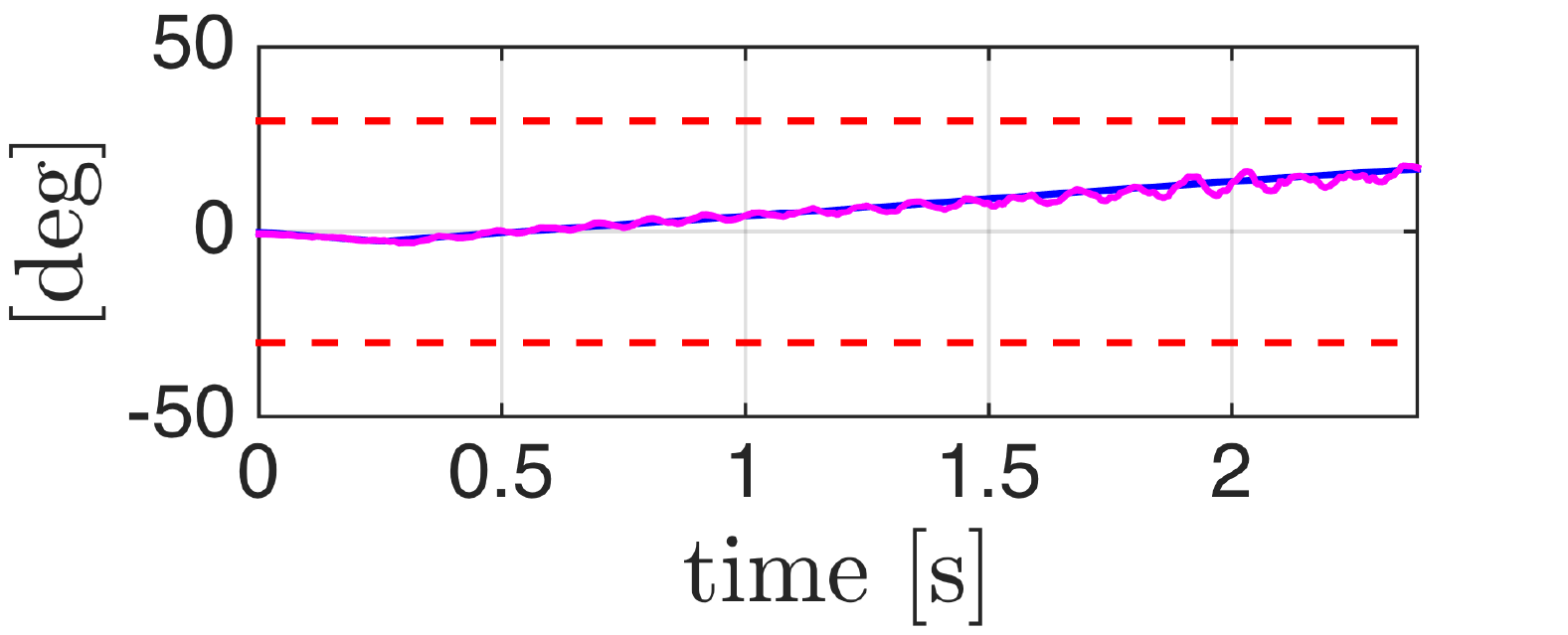} \label{fig:phi_exp}}}}\\
		&
		\hspace{-0.7cm}{\subfloat[][Velocity norm $||\text{\textbf{v}}||$]{{\includegraphics[width=4.7cm]{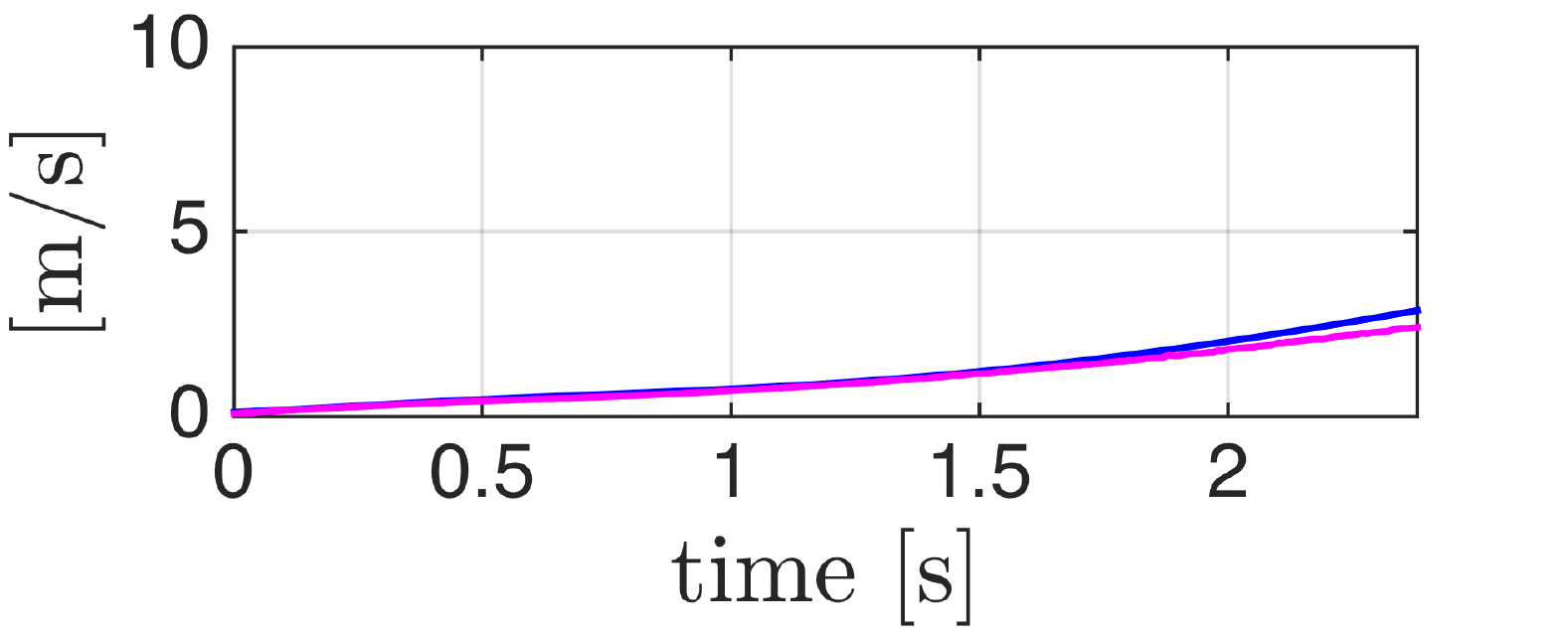} \label{fig:vel_exp}}}}\\
		\hspace{0.3cm}{\subfloat[][Angular rate $p$]{\includegraphics[width=4.7cm]{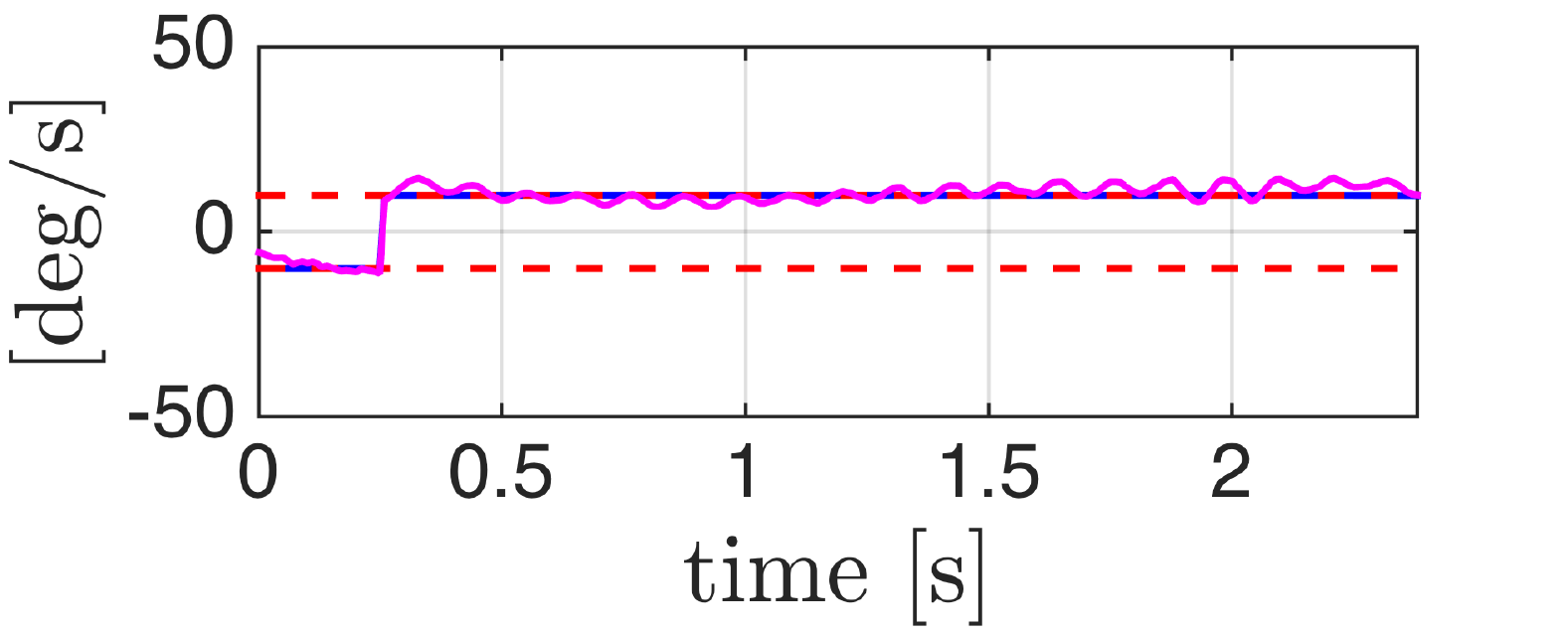}\vspace{0.2cm}\label{fig:pp_exp}}} &
		\hspace{-0.4cm}{\subfloat[][Thrust $\force$]{\hspace{-0.2cm}\includegraphics[width=4.7cm]{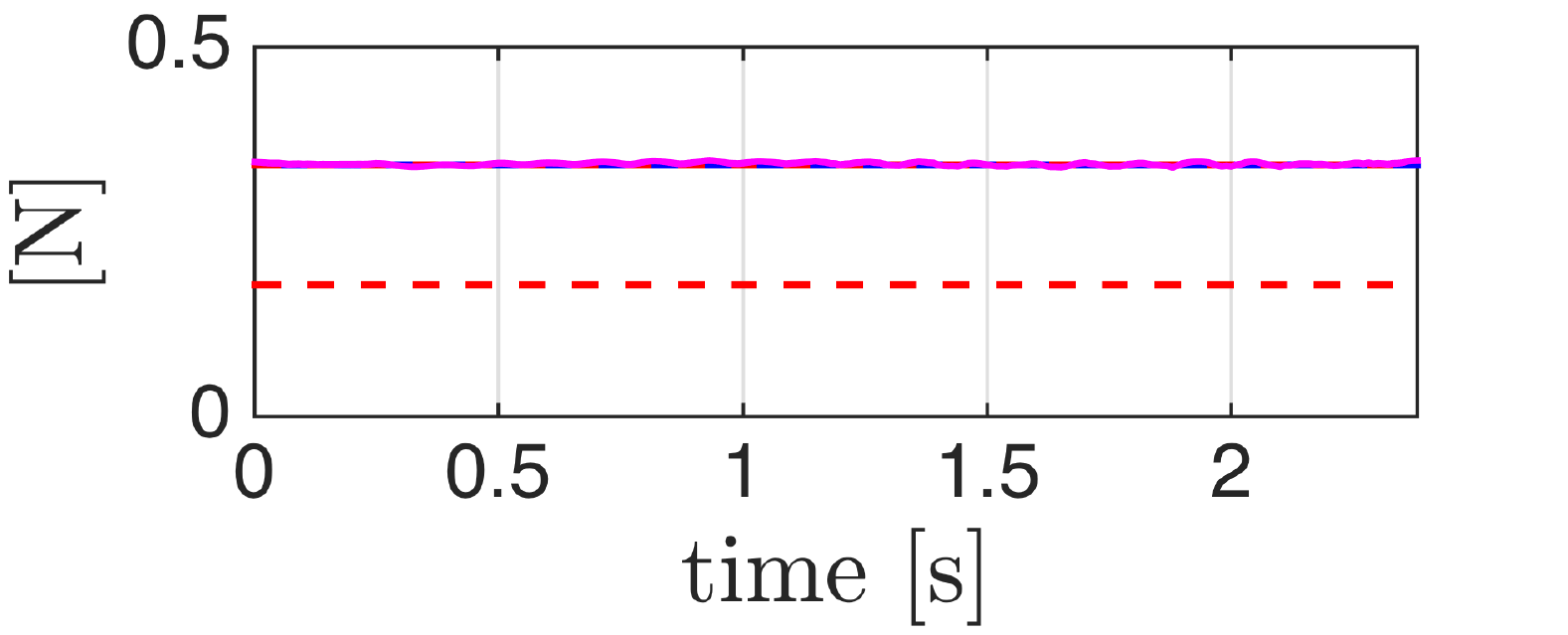}\vspace{0.2cm}\label{fig:ff_exp}}}\\
		\multicolumn{2}{c}{\subfloat[][Experiment snapshots]{\includegraphics[scale=0.15]{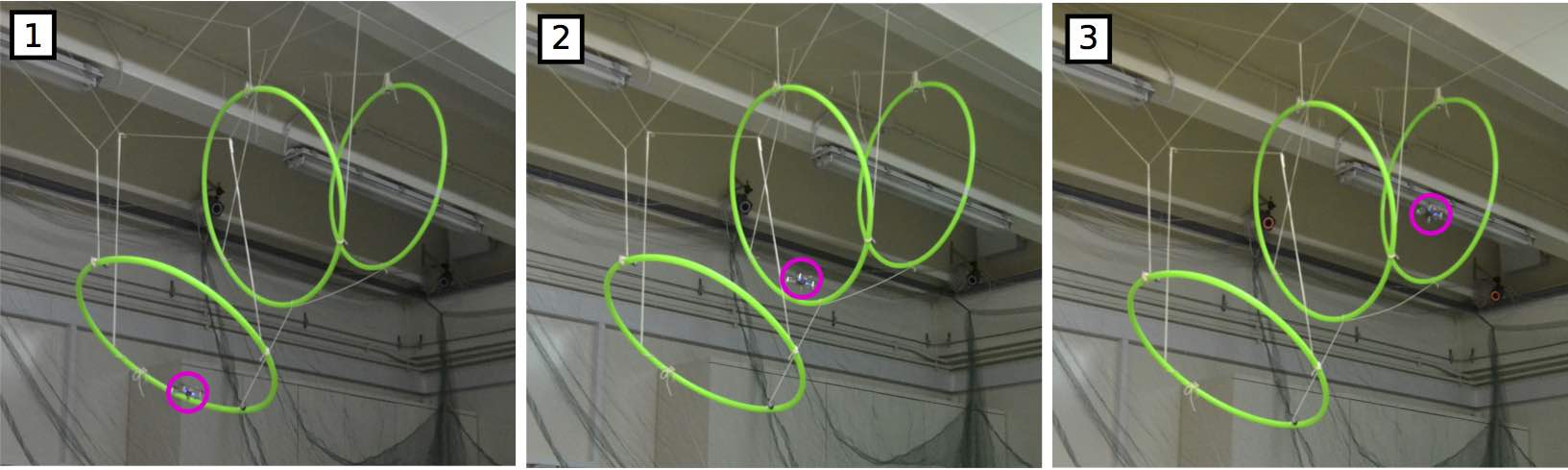}\label{fig:exp_snapshots}}}
	\end{tabular}
	\caption{Experimental test. Desired trajectory (blue) and actual trajectory (magenta). Constraint boundaries are depicted in dashed red. Hula hoops are depicted in solid green.}
	\label{fig:exp}
\end{figure}

{By using our minimum-time strategy, we obtain the following result.}
The optimal trajectory, performed in $2.38$ s,
is depicted in Figure \ref{fig:exp} (solid blue). Constraint boundaries are depicted in dashed red and the hula hoops are depicted in solid green. 
The optimal path (blue line in Figure \ref{fig:path_exp}) first takes negative values of $p_2$ until changing direction toward positive $p_2$ values, touching the constraint boundary in the proximity of the maximum curvature of the tube, and staying in the middle of the feasibility region at the end.
The roll angle (blue line in Figure \ref{fig:phi_exp}) decreases in order to push the vehicle to negative $p_2$ values and then it monotonically increases during the remaining time interval.
The velocity module (blue line in Figure \ref{fig:vel_exp}) always increases, as we expect for a minimum-time trajectory.
As regards the inputs, in the beginning, the angular rate $p$ (blue line in Figure \ref{fig:pp_exp}) stays on the lower bound and then it switches to the upper bound.
The thrust $\force$ (blue line in Figure \ref{fig:ff_exp}) always takes the upper bound.

We execute the computed minimum-time trajectory on the CrazyFlie nano-quadrotor
by using the closed-loop, maneuver regulation controller developed in
\cite{SS-GN-HHB-AF:13}, in which the minimum-time trajectory is used for the
desired maneuver.
The maneuver regulation controller computes thrust and angular rate virtual
inputs, which are tracked by the standard off-the-shelf angular rate
controller provided on board the CrazyFlie.
The actual (experimental) trajectory performed using our maneuver regulation
controller is depicted in Figure \ref{fig:exp} in solid magenta. Snapshots of
the experiment are depicted in Figure \ref{fig:exp_snapshots}.  As expected,
the quadrotor passes close to the second hula hoop maintaining the distance
imposed by the restrictive constraints in the optimization problem. 
The actual velocity does not perfectly match (at higher velocities) the desired one, due to the unmodeled drag effect.
Since the vehicle is asked to follow the desired thrust, the actual velocity becomes lower than the desired one because of the opposing aerodynamic force.
The experiment shows the actual feasibility of the optimal trajectory and also reveals that a more accurate model including aerodynamic effects would improve the control performance.

\section{Conclusion}
In this paper, we have presented a strategy to address the minimum-time problem
for quadrotors in constrained environments. Our approach consists of: (i)
generating a frame path, (ii) expressing the quadrotor dynamics in a new
set of coordinates ``transverse" with respect to that path, and (iii) redefining
cost and constraints in the new coordinates.  Thus, we obtain a
reformulation of the problem, which we solve by combining the PRONTO algorithm
with a barrier function approach.  Numerical computations on two challenging
scenarios prove the effectiveness of the strategy and allow us to show
interesting dynamic capabilities of the vehicle. Moreover, the experimental test
of the second scenario shows the feasibility of the computed
trajectory. As a future work, we aim at extending our strategy to a
  scenario with moving obstacles. Challenges to be addressed include how to
  combine trajectory generation and control, and how to take into account a fast
  integration of the dynamics for realtime computation.

\appendices

\section{Projection Operator Newton Method}
\label{app:pronto}
Here, we provide a brief description of the PRONTO algorithm \cite{hauser2002projection}.
The PRONTO algorithm is based on a properly designed \emph{projection operator}
$\mathcal{P} : \xi_c \rightarrow \xi$, mapping a state-control curve
$\xi_c =(\bar{\boldsymbol{x}}_{w,c}(\cdot), \bar{\boldsymbol{u}}_c(\cdot))$ into a system
trajectory 
$\xi =(\bar{\boldsymbol{x}}_w(\cdot), \bar{\boldsymbol{u}}(\cdot))$, by the nonlinear feedback system
\begin{align}
\bar{\boldsymbol{x}}'_w(s) & = \bar{f}(\bar{\boldsymbol{x}}_w(s),\bar{\boldsymbol{u}}(s)), \quad \bar{\boldsymbol{x}}_w(0) = \boldsymbol{x}_{w0}, \nonumber\\
  \bar{\boldsymbol{u}}(s) & = \bar{\boldsymbol{u}}_c(s) + \bar{K}(s)(\bar{\boldsymbol{x}}_{w,c}(s)-\bar{\boldsymbol{x}}_w(s)),
                     \label{eq:proj_oper_def}
\end{align}
where the feedback gain $\bar{K}(\cdot)$ is designed by solving a suitable
linear quadratic optimal control problem on the linearized dynamics of \eqref{eq:tran_dynamics} about the trajectory $\xi$. 
Note that the feedback gain $\bar{K}(\cdot)$ is only used to define the projection operator and it is not related to the controller used to execute the optimal trajectory in our experimental test.
The projection operator is used to
convert the dynamically constrained optimization problem \eqref{eq:mintime2}
into the unconstrained problem
\begin{equation}
  \begin{split}
    \min_{\xi} &\; g(\xi;\bar{k}),
  \end{split}
  \label{eq:mintime3}
\end{equation}
where
$g(\xi;\bar{k}) = h(\mathcal{P}(\xi); \bar{k})$, and $h(\xi; \bar{k}) := \int_0^L \!
(\frac{1-\bar{k}(s) \bar{w}_1(s)}{\bar{\boldsymbol{t}}^T(s) \bar{\text{\textbf{v}}}(s)} + \epsilon \sum_j \beta_\nu (-c_j(\bar{\boldsymbol{x}}_w(s),\bar{\boldsymbol{u}}(s)))) ds +
\; \epsilon_f \sum_i \beta_{\nu_f} (-c_{f,i}(\bar{\boldsymbol{x}}_w(L)))$.
Then, using an (infinite dimensional) Newton descent method, a local minimizer
of \eqref{eq:mintime3} is computed iteratively. 
Given the current trajectory iterate $\xi_i$, the search
direction $\zeta_i$ is obtained by solving a linear quadratic optimal control
problem with cost
$Dg(\xi_i; \bar{k}) \cdot \zeta + \frac{1}{2} D^2 g(\xi_i;
\bar{k})(\zeta,\zeta)$,
where $\zeta \mapsto Dg(\xi_i; \bar{k}) \cdot \zeta$ and
$\zeta \mapsto D^2 g(\xi_i; \bar{k})(\zeta,\zeta)$ are respectively the first
and second Fr\'echet differentials of the functional $g(\xi,\bar{k})$ at
$\xi_i$. Then, the curve $\xi_i + \gamma_i \zeta_i$, where $\gamma_i$ is a step
size obtained through a standard backtracking line search, is projected, by
means of the projection operator, in order to get a new trajectory $\xi_{i+1}$.

The strength of this approach is that the local minimizer of \eqref{eq:mintime3}
is obtained as the limit of a sequence of trajectories, i.e., curves satisfying
the dynamics.  Furthermore, the feedback system \eqref{eq:proj_oper_def},
defining the projection operator, allows us to generate trajectories in a
numerically stable manner.

\begin{remark}
An elegant extension of the PRONTO method to Lie groups is developed in  \cite{saccon2013optimal} and could be alternatively used in our strategy.
\end{remark}

\ifCLASSOPTIONcaptionsoff
  \newpage
\fi

\bibliographystyle{IEEEtran}
\bibliography{bibliography}  

\begin{thebibliography}{10}
\providecommand{\url}[1]{#1}
\csname url@samestyle\endcsname
\providecommand{\newblock}{\relax}
\providecommand{\bibinfo}[2]{#2}
\providecommand{\BIBentrySTDinterwordspacing}{\spaceskip=0pt\relax}
\providecommand{\BIBentryALTinterwordstretchfactor}{4}
\providecommand{\BIBentryALTinterwordspacing}{\spaceskip=\fontdimen2\font plus
\BIBentryALTinterwordstretchfactor\fontdimen3\font minus
  \fontdimen4\font\relax}
\providecommand{\BIBforeignlanguage}[2]{{%
\expandafter\ifx\csname l@#1\endcsname\relax
\typeout{** WARNING: IEEEtran.bst: No hyphenation pattern has been}%
\typeout{** loaded for the language `#1'. Using the pattern for}%
\typeout{** the default language instead.}%
\else
\language=\csname l@#1\endcsname
\fi
#2}}
\providecommand{\BIBdecl}{\relax}
\BIBdecl

\bibitem{dadkhah2012survey}
N.~Dadkhah and B.~Mettler, ``{S}urvey of motion planning literature in the
  presence of uncertainty: considerations for {UAV} guidance,'' \emph{Journal
  of Intelligent \& Robotic Systems}, vol.~65, no. 1-4, pp. 233--246, 2012.

\bibitem{bottasso2008path}
C.~L. Bottasso, D.~Leonello, and B.~Savini, ``{P}ath {P}lanning for
  {A}utonomous {V}ehicles by {T}rajectory {S}moothing {U}sing {M}otion
  {P}rimitives,'' \emph{IEEE Transactions on Control Systems Technology},
  vol.~16, no.~6, pp. 1152--1168, 2008.

\bibitem{ambrosino2009path}
G.~Ambrosino, M.~Ariola, U.~Ciniglio, F.~Corraro, E.~D. Lellis, and A.~Pironti,
  ``{P}ath {G}eneration and {T}racking in 3-{D} for {UAV}s,'' \emph{IEEE
  Transactions on Control Systems Technology}, vol.~17, no.~4, pp. 980--988,
  2009.

\bibitem{herisse2012landing}
B.~Heriss{\'e}, T.~Hamel, R.~Mahony, and F.-X. Russotto, ``{L}anding a {VTOL}
  unmanned aerial vehicle on a moving platform using optical flow,'' \emph{IEEE
  Transactions on Robotics}, vol.~28, no.~1, pp. 77--89, 2012.

\bibitem{naldi2015robust}
R.~Naldi, A.~Torre, and L.~Marconi, ``Robust control of a miniature ducted-fan
  aerial robot for blind navigation in unknown populated environments,''
  \emph{IEEE Transactions on Control Systems Technology}, vol.~23, no.~1, pp.
  64--79, 2015.

\bibitem{hou2016dynamic}
X.~Hou and R.~Mahony, ``{D}ynamic {K}inesthetic {B}oundary for {H}aptic
  {T}eleoperation of {VTOL} {A}erial {R}obots in {C}omplex {E}nvironments,''
  \emph{IEEE Transactions on Systems, Man, and Cybernetics: Systems}, vol.~46,
  no.~5, pp. 694--705, 2016.

\bibitem{furci2015plan}
M.~Furci, R.~Naldi, S.~Karaman, and L.~Marconi, ``{A} {C}ombined {P}lanning and
  {C}ontrol {S}trategy for {M}obile {R}obots {N}avigation in {P}opulated
  {E}nvironments,'' in \emph{IEEE Conference on Decision and Control}, 2015.

\bibitem{masoud2015plan}
A.~Masoud and A.~Al-Shaikhi, ``{T}ime-sensitive, sensor-based, joint planning
  and control of mobile robots in cluttered spaces: {A} harmonic potential
  approach,'' in \emph{IEEE Conference on Decision and Control}, 2015.

\bibitem{cowling2010direct}
I.~D. Cowling, O.~A. Yakimenko, J.~F. Whidborne, and A.~K. Cooke, ``{D}irect
  method based control system for an autonomous quadrotor,'' \emph{Journal of
  Intelligent \& Robotic Systems}, vol.~60, no.~2, pp. 285--316, 2010.

\bibitem{mellinger2011minimum}
D.~Mellinger and V.~Kumar, ``{M}inimum snap trajectory generation and control
  for quadrotors,'' in \emph{IEEE International Conference on Robotics and
  Automation}, 2011.

\bibitem{bouktir2008trajectory}
Y.~Bouktir, M.~Haddad, and T.~Chettibi, ``{T}rajectory planning for a quadrotor
  helicopter,'' in \emph{Mediterranean Conference on Control and Automation},
  2008.

\bibitem{van2013time}
W.~Van~Loock, G.~Pipeleers, and J.~Swevers, ``{T}ime-optimal quadrotor
  flight,'' in \emph{European Control Conference}, 2013.

\bibitem{chen2016online}
J.~Chen, T.~Liu, and S.~Shen, ``Online generation of collision-free
  trajectories for quadrotor flight in unknown cluttered environments,'' in
  \emph{IEEE International Conference on Robotics and Automation}, 2016.

\bibitem{bry2015aggressive}
A.~Bry, C.~Richter, A.~Bachrach, and N.~Roy, ``{A}ggressive flight of
  fixed-wing and quadrotor aircraft in dense indoor environments,'' \emph{The
  International Journal of Robotics Research}, vol.~34, no.~7, pp. 969--1002,
  2015.

\bibitem{koyuncu2008probabilistic}
E.~Koyuncu and G.~Inalhan, ``{A} probabilistic b-spline motion planning
  algorithm for unmanned helicopters flying in dense 3d environments,'' in
  \emph{IEEE/RSJ International Conference on Intelligent Robots and Systems},
  2008.

\bibitem{bouffard2009hybrid}
P.~M. Bouffard and S.~L. Waslander, ``A hybrid randomized/nonlinear programming
  technique for small aerial vehicle trajectory planning in 3d,''
  \emph{Planning, Perception and Navigation for Intelligent Vehicles}, vol.~63,
  2009.

\bibitem{devaurs2015optimal}
D.~Devaurs, T.~Sim{\'e}on, and J.~Cort{\'e}s, ``{O}ptimal {P}ath {P}lanning in
  {C}omplex {C}ost {S}paces {W}ith {S}ampling-{B}ased {A}lgorithms,''
  \emph{IEEE Transactions on Automation Science and Engineering}, vol.~13,
  no.~2, pp. 415--424, 2016.

\bibitem{allen2016real}
R.~Allen and M.~Pavone, ``{A} {R}eal-{T}ime {F}ramework for {K}inodynamic
  {P}lanning with {A}pplication to {Q}uadrotor {O}bstacle {A}voidance,'' in
  \emph{AIAA Conf. on Guidance, Navigation and Control, San Diego, CA}, 2016.

\bibitem{hehn2015real}
M.~Hehn and R.~D'Andrea, ``{R}eal-{T}ime {T}rajectory {G}eneration for
  {Q}uadrocopters,'' \emph{IEEE Transactions on Robotics}, vol.~31, no.~4, pp.
  877--892, 2015.

\bibitem{augugliaro2012generation}
F.~Augugliaro, A.~P. Schoellig, and R.~D'Andrea, ``{G}eneration of
  collision-free trajectories for a quadrocopter fleet: {A} sequential convex
  programming approach,'' in \emph{IEEE/RSJ International Conference on
  Intelligent Robots and Systems}, 2012.

\bibitem{augugliaro2013dance}
------, ``{D}ance of the flying machines: {M}ethods for designing and executing
  an aerial dance choreography,'' \emph{IEEE Robotics \& Automation Magazine},
  vol.~20, no.~4, pp. 96--104, 2013.

\bibitem{hauser2002projection}
J.~Hauser, ``{A} projection operator approach to the optimization of trajectory
  functionals,'' in \emph{IFAC world congress}, 2002.

\bibitem{hauser2006barrier}
J.~Hauser and A.~Saccon, ``{A} barrier function method for the optimization of
  trajectory functionals with constraints,'' in \emph{IEEE Conference on
  Decision and Control}, 2006.

\bibitem{hauser2006motorcycle}
------, ``Motorcycle modeling for high-performance maneuvering,'' \emph{IEEE
  Control Systems}, vol.~26, no.~5, pp. 89--105, 2006.

\bibitem{hauslerenergy}
A.~J. Hausler, A.~Saccon, A.~P. Aguiar, J.~Hauser, and A.~M. Pascoal,
  ``{E}nergy-{O}ptimal {M}otion {P}lanning for {M}ultiple {R}obotic {V}ehicles
  {W}ith {C}ollision {A}voidance,'' \emph{IEEE Transactions on Control Systems
  Technology}, vol.~24, no.~3, pp. 867--883, 2015.

\bibitem{rucco2015virtual}
A.~Rucco, A.~P. Aguiar, and J.~Hauser, ``{A} {V}irtual {T}arget {A}pproach for
  {T}rajectory {O}ptimization of a {G}eneral {C}lass of {C}onstrained
  {V}ehicles,'' in \emph{IEEE Conference on Decision and Control}, 2015.

\bibitem{hua2013introduction}
M.~D. Hua, T.~Hamel, P.~Morin, and C.~Samson, ``{I}ntroduction to feedback
  control of underactuated {VTOL} vehicles: {A} review of basic control design
  ideas and principles,'' \emph{IEEE Control Systems}, vol.~33, no.~1, pp.
  61--75, 2013.

\bibitem{PSM:08}
M.~P. Setterlund, ``{G}eometric-based {S}patial {P}ath {P}lanning,'' PhD
  dissertation, University of Texas at Austin, 2008.

\bibitem{AB-JH:94}
A.~Banaszuk and J.~Hauser, ``Feedback linearization of transverse dynamics for
  periodic orbits,'' in \emph{IEEE Conference on Decision and Control}, 1994.

\bibitem{AS-JH-AB:13}
A.~Saccon, J.~Hauser, and A.~Beghi, ``{A} {V}irtual {R}ider for {M}otorcycles:
  {M}aneuver {R}egulation of a {M}ulti-{B}ody {V}ehicle {M}odel,'' \emph{IEEE
  Transactions on Control Systems Technology}, vol.~21, no.~2, pp. 332--346,
  March 2013.

\bibitem{SS-GN-HHB-AF:13}
S.~Spedicato, G.~Notarstefano, H.~H. B{\"u}lthoff, and A.~Franchi, ``Aggressive
  {M}aneuver {R}egulation of a {Q}uadrotor {UAV},'' in \emph{The 16th
  International Symposium on Robotics Research}, 2013.

\bibitem{saccon2013optimal}
A.~Saccon, J.~Hauser, and A.~P. Aguiar, ``{O}ptimal control on {L}ie groups:
  {T}he projection operator approach,'' \emph{IEEE Transactions on Automatic
  Control}, vol.~58, no.~9, pp. 2230--2245, 2013.

\end{thebibliography}

\begin{IEEEbiography}[{\includegraphics[width=1in]{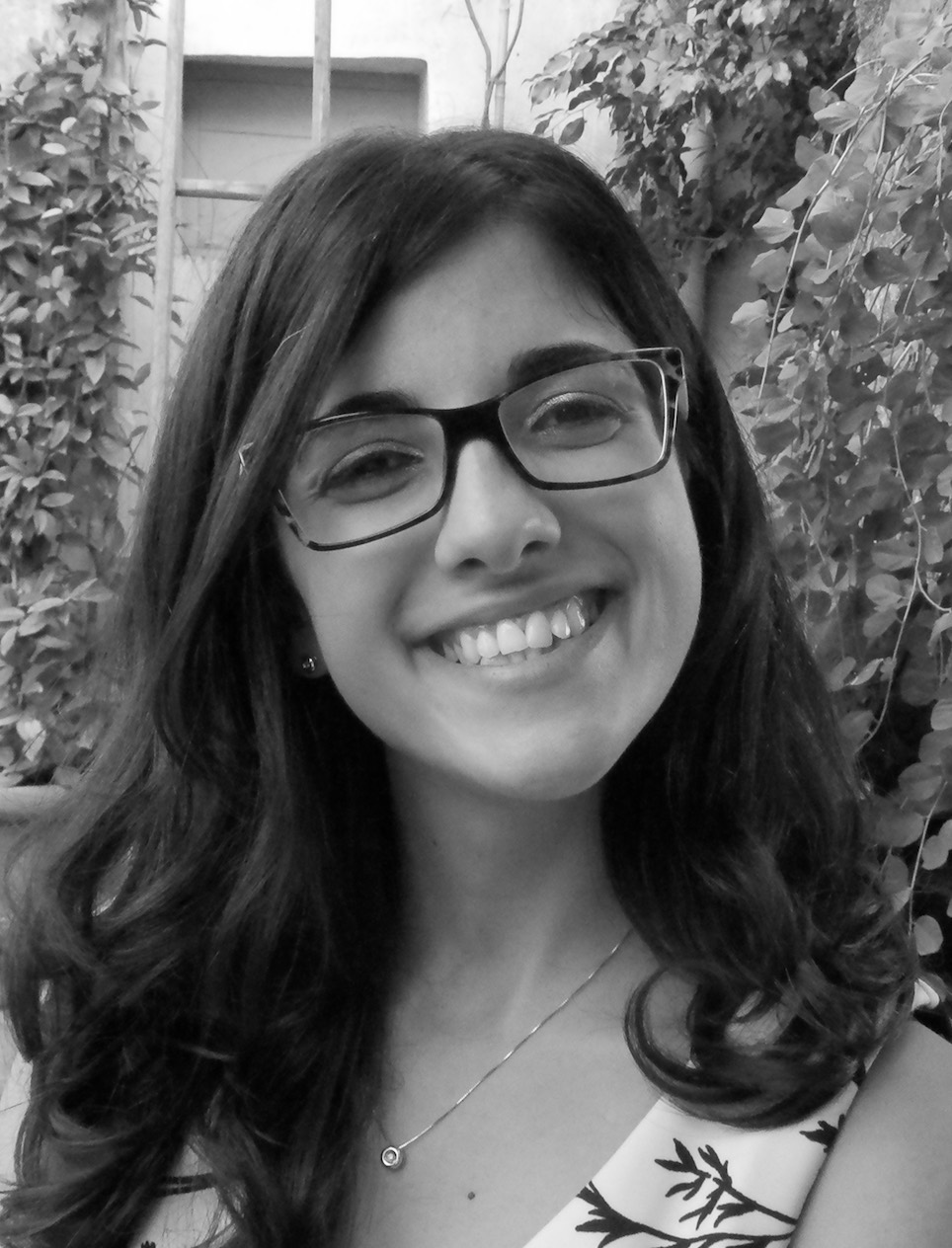}}]
{Sara Spedicato}
is a Post-Doctoral Researcher at the Universit\`a del Salento (Lecce, Italy), where she 
received the Laurea degree ``summa cum laude" in Mechanical Engineering (curriculum ``Servomechanisms and industrial automation") in 2012 and the Ph.D. degree in Information Engineering in 2016.
She carried out an internship activity at the ETH Z\"urich (Z\"urich, Switzerland) from June to September 2012. She was a visiting graduate student at the Max Planck Institute for Biological Cybernetics (T\"ubingen, Germany) from January to August 2013. Her research activity involves nonlinear optimal control, distributed optimization, trajectory optimization and maneuvering for autonomous aerial vehicles.
\end{IEEEbiography}

\begin{IEEEbiography}[{\includegraphics[width=1in]{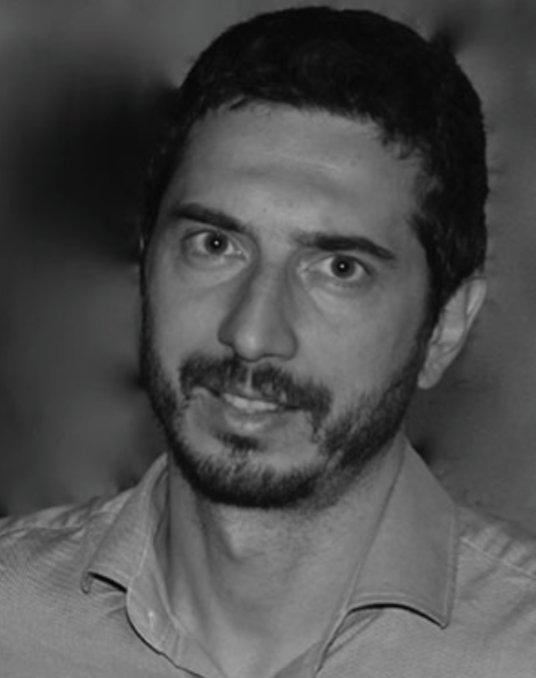}}]
{Giuseppe Notarstefano}(M'11) is Associate Professor at the Universit\`a del Salento (Lecce, Italy), where he was Assistant Professor (Ricercatore) from February 2007 to May 2016. He received the Laurea degree ``summa cum laude" in Electronics Engineering from the Universit\`a di Pisa in 2003 and the Ph.D. degree in Automation and Operation Research from the Universit\`a di Padova in 2007. He has been visiting scholar at the University of Stuttgart, University of California Santa Barbara and University of Colorado Boulder. His research interests include distributed optimization, cooperative control in complex networks, applied nonlinear optimal control, and trajectory optimization and maneuvering of aerial and car vehicles.
He serves as an Associate Editor for IEEE Transactions on Control Systems Technology, for IEEE Control Systems Letters, for the Conference Editorial Board of the IEEE Control Systems Society and for other IEEE and IFAC conferences.
He coordinated the VI-RTUS team winning the International Student Competition Virtual Formula 2012. He is recipient of an ERC Starting Grant 2014.
\end{IEEEbiography}

\end{document}